\input{amssymb.sty}
\documentclass[11pt]{amsart}
\theoremstyle{plain}
\newtheorem{thm}{Theorem}
\newtheorem*{thm*}{Theorem}
\newtheorem{prop}{Proposition}
\newtheorem{lemma}{Lemma}
\newtheorem{cor}{Corollary}
\newtheorem*{conj*}{Conjecture}

\theoremstyle{definition}
\newtheorem*{defn*}{Definition}
\newtheorem{rems}[thm]{Remarks}
\newtheorem*{rems*}{Remarks}
\newtheorem*{proof*}{Proof}
\newcommand{\npartial}{\slash\!\!\!\partial}
\newcommand{\Heis}{\operatorname{Heis}}
\newcommand{\Solv}{\operatorname{Solv}}
\newcommand{\Spin}{\operatorname{Spin}}

\newcommand{\ind}{\operatorname{ind}}
\newcommand{\Index}{\operatorname{index}}
\newcommand{\ch}{\operatorname{ch}}
\newcommand{\rank}{\operatorname{rank}}

\newcommand{\ip}[2]{\langle #1,#2\rangle}
\newcommand{\nc}{\newcommand}
\nc{\nt}{\newtheorem}
\nc{\gf}[2]{\genfrac{}{}{0pt}{}{#1}{#2}}
\nc{\mb}[1]{{\mbox{$ #1 $}}}
\nc{\real}{{\mathbb R}}
\nc{\comp}{{\mathbb C}}
\nc{\ints}{{\mathbb Z}}
\nc{\Ltoo}{\mb{L^2({\mathbf H})}}
\nc{\rtoo}{\mb{{\mathbf R}^2}}
\nc{\slr}{{\mathbf {SL}}(2,\real)}
\nc{\slz}{{\mathbf {SL}}(2,\ints)}
\nc{\su}{{\mathbf {SU}}(1,1)}
\nc{\so}{{\mathbf {SO}}}
\nc{\hyp}{{\mathbb H}}
\nc{\disc}{{\mathbf D}}
\nc{\torus}{{\mathbb T}}

\newcommand{\tr}{\mbox{tr}}

\nc{\ca}{{\mathcal A}}
\nc{\cag}{{{\mathcal A}^\Gamma}}
\nc{\cg}{{\mathcal G}}
\nc{\chh}{{\mathcal H}}
\nc{\ck}{{\mathcal B}}
\nc{\cl}{{\mathcal L}}
\nc{\cm}{{\mathcal M}}
\nc{\cs}{{\mathcal S}}
\nc{\cz}{{\mathcal Z}}
\nc{\sind}{\sigma{\rm -ind}}
\begin{document}

\title[Quantum Hall Effect on the Hyperbolic Plane]
{Quantum Hall Effect on the Hyperbolic Plane}
\author{A. L. Carey}
\address{Department of Mathematics, University of Adelaide, Adelaide 5005,
Australia}
\email{acarey@maths.adelaide.edu.au}
\author{K. C. Hannabuss}
\address{Department of Mathematics, University of Oxford, England.}
\email{khannabu@maths.adelaide.edu.au}
\author{V. Mathai}
\address{Department of Mathematics, University of Adelaide, Adelaide 5005,
Australia}
\email{vmathai@maths.adelaide.edu.au}
\author{P. McCann}
\address{Department of Mathematics, University of Adelaide, Adelaide 5005,
Australia}
\email{pmccann@maths.adelaide.edu.au}

\date{FEBRUARY 1997}

\begin{abstract}
In this paper, we study both the continuous model and the discrete model of the
Quantum Hall Effect (QHE) on the hyperbolic plane.
The Hall conductivity is identified as a geometric invariant associated to an
imprimitivity algebra of observables. We define a twisted analogue of the
Kasparov map, which
enables us to use the pairing between $K$-theory and cyclic cohomology theory,
to identify this geometric invariant with a topological index, thereby
proving the integrality of the Hall conductivity in this case.
\end{abstract}

\maketitle


\section{Introduction}

The usual model of the integer quantum Hall effect involves electrons moving
in a two dimensional conductor under the influence of a magnetic field.
The field is applied in a direction orthogonal to the
conductor. The Hamiltonian used is that for a single electron moving
under the influence of
 this magnetic field with the addition of a potential term to represent
the field due to the
lattice of ions making up the conductor. For simplicity this
lattice is often assumed to be
periodic in the two axis directions in the plane.
The effect of impurities can be modelled by departing from a perfectly
periodic potential. The definitive treatment from a mathematical point of
 view is due to Bellisard \cite{Bel+E+S} and Xia \cite{Xia}. In this approach
no assumption is
made about the rationality of the imposed magnetic flux while the integrality
 of the Hall conductance follows by showing that it is given by
the index of a Fredholm operator. Xia  also exhibits the conductance
as a  topological index.
These demonstrations use in an
essential way Connes' non-commutative differential geometry.

In this paper we are interested in what can be said when
one replaces the usual two dimensional conducting material with its
Euclidean geometry by a two dimensional sample with hyperbolic geometry.
Physically one should think of hyperbolic space
and hence the sample as an embedded hyperboloid
in Euclidean 3-space. The crystal lattice of the conductor is now modelled
by the orbit of a freely acting discrete group.
 For reasons of convenience we take
this to be the fundamental group of a Riemann surface (though aspects of
our analysis work more generally). The magnetic field remains orthogonal to
the two dimensional conductor and the electric potential we take to be
periodic under the action of this discrete group (we do not attempt to model
impurities). When the magnetic flux is rational, spectral
properties of Hamiltonians for
a single electron have been studied, particularly in the cases where
the Hamiltonian may be defined on a finite cover of the Riemann
surface \cite{Iengo+Li}. In this paper however
 we are interested in approaching the problem using non-commutative
 geometry, which allows the flux to be any
real number, and the algebras of interest are not always associated
with vector bundles over the Riemann surface.

We begin by reviewing the construction of the Hamiltonian.
First we take as our principal model of hyperbolic space,
the hyperbolic plane.
This is the upper half-plane $\hyp$ in $\comp$ equipped with its usual
Poincar\'e metric $(dx^2+dy^2)/y^2$, and symplectic area form
$\omega_\hyp  = dx\wedge dy/y^2$.
The group $\slr$ acts transitively on $\hyp$ by M\"obius transformations
$$
x+iy = \zeta  \mapsto g\zeta = \frac{a\zeta+b}{c\zeta+d},\quad\mbox{for }
g=\left(
\begin{array}{cc}
a & b\\
c & d
\end{array}\right).$$
Any Riemann surface of genus $g$ greater than 1 can be realised as the
quotient of
$\hyp$ by the action of its fundamental group realised as a subgroup
$\Gamma$ of $\slz \subset \slr$.

Let us now pick a 1-form $\eta$ such that $d\eta = \theta\omega_\hyp$, for some
fixed $\theta \in (0,1]$.
As in geometric quantisation we may regard $\eta$ as defining a connection
$\nabla = d-i\eta$ on a line bundle $\cl$ over $\hyp$, whose curvature is
$\theta\omega_\hyp$.
Physically we can think of $\eta$ as the electromagnetic vector potential for
 a uniform magnetic field of strength $\theta$ normal to $\hyp$.
Using the Riemannian metric  the Hamiltonian of an
electron in this field is given in suitable units by
$$H = H_\eta = -\frac 12\nabla^*\nabla = -\frac 12(d-i\eta)^*(d-i\eta).$$
Comtet [Comtet] has shown that $H$ differs from a multiple of the Casimir
element for $\slr$, $\frac 18{\bf J}.{\bf J}$,
by a constant, where $J_1$, $J_2$ and $J_3$
denote a certain
representation of generators of the Lie algebra $sl(2,\real)$, satisfying
$$[J_1,J_2] = -iJ_3, \qquad [J_2,J_3] = iJ_1, \qquad [J_3,J_1] = iJ_2,$$
so that ${\bf J}.{\bf J} = J_1^2+J_2^2-J_3^2$ is the quadratic Casimir element.
This shows very clearly the underlying $\slr$-invariance of the theory.
In a real material this Hamiltonian would be modified by the addition of a
potential $V$.
By taking $V$ to be invariant under $\Gamma$ this perturbation is
given a crystalline type structure analogous to the use of periodic potentials
invariant under $\ints^2$ in the Euclidean plane $\real^2$.
Comtet has computed the spectrum of the unperturbed Hamiltonian $H_\eta$
for $\eta = -\theta dx/y$ to be the union of finitely many eigenvalues
$\{(2k+1)\theta -k(k+1):k=0,1,2\ldots < \theta-\frac 12\}$,
and the continuous spectrum
$[\frac 14 + \theta^2, \infty)$ for those values of $\theta$
for which the de Rham
cohomology class of $\theta\omega_\hyp$ is integral. Its zeta
function and the kernel of its resolvent are also known in this case,
\cite{Comtet}, \cite{Comtet+H}.
Any $\eta$ is cohomologous to $-\theta dx/y$ (since they both have
$\omega_\hyp$ as differential) and forms differing by an exact form $d\phi$
 give
equivalent models: in fact, multiplying the wave functions by
$\exp(i\phi)$ shows that the models for $\eta$ and $-\theta dx/y$ are
unitarily equivalent.
This equivalence also intertwines the $\Gamma$-actions so that the spectral
densities for the two models also coincide.
However, it is the perturbed Hamiltonian $H_{\eta,V} = H_\eta +V$ which
is the key to the quantum Hall effect on the hyperbolic plane, and the spectrum
of this is unknown for general $\Gamma$-invariant $V$. As we noted above
the Hall effect on Riemann surfaces has also been considered \cite{Iengo+Li},
\cite{Av+K+P+S} but this is different from the problem we consider here.

These considerations suggest that one could mimic the
non-commutative geometry approach of Bellissard-Connes
to the integer quantum Hall effect on Euclidean space
\cite{Bel}, \cite{Nak+Bel}, \cite{Bel+E+S}, \cite{Co2}, \cite{Xia}
 in a hyperbolic setting.
This interprets the Hall conductivity as a non-commutative Chern character,
whose integrality follows from K-theory.
Physically such situations have been considered without the perturbing
potential (or $\Gamma$ is trivial) in the context of exploring edge effects for
the quantum Hall effect and the behaviour of electrons in quantum dots.
Much of the mathematical machinery needed for this has already been discussed
in a geometrical context \cite{Co}, \cite{Co2}, \cite{Comtet},
\cite{Comtet+H} and will be exploited
here.
We  also discuss the discrete version of the theory \cite{Co2},
\cite{Sun} motivated in part by
some results of \cite{MC} in the Euclidean
setting. These hyperbolic Hall effect models occupy the first seven sections.
Specifically we show that there is a principal groupoid $C^*$-algebra
with cocycle for the diagonal $\Gamma$ action on  $\hyp\times\hyp$
which contains the resolvent of the various  Hamiltonians we consider.
In order to construct a Fredholm module for this algebra we found it
useful to take a more abstract group theoretic approach. We show that
our groupoid algebra is isomorphic to a quotient of the
$\Gamma$ invariant part of the imprimitivity algebra for inducing
from the maximal compact subgroup of $\slr$ to $C^*(\slr,\sigma)$
(the multiplier, or group 2-cocycle, $\sigma$ extends to all of $\slr$). This
imprimitivity algebra has a regular representation, induced by a
canonical trace, the Hilbert space of which provides a Fredholm module
which is 2-summable for a dense subalgebra of the imprimitivity algebra.
We show that this dense subalgebra contains the spectral projections
corresponding to gaps in the spectrum of our Hamiltonians. Similar
results hold in the discrete model as well.
The  connection between the continuous and discrete models
arises from the Morita
 equivalence of our quotient of the $\Gamma$ invariant
imprimitivity algebra with $C^*(\Gamma,\sigma)$.

The main results of our paper follow by extending the approach of \cite{Xia}
to cover the hyperbolic case. In fact in Section 8
we prove some general theorems about the $K$-groups of
$C^*$-algebras which generalise those arising from the hyperbolic Hall effect.
The relevance of $K$-theory
can be understood in the case of the integer Hall effect
on Euclidean space partly as a result of the calculation \cite{Elliott},
 \cite{Bel}, \cite{Co}:
$$
K_*(C^*(\mathbb Z^n,\sigma)) \cong K_*(C^*(\mathbb Z^n)) \cong K^*(\mathbb T^n)
$$
for any multiplier
({\it i.e.} group 2-cocycle) $\sigma$ on $\mathbb Z^n$.
This result
 has lead to the twisted group $C^*$-algebras
$C^*(\mathbb Z^n,\sigma)$ being called {\em noncommutative tori}.
 This calculation was generalized by
Packer and Raeburn \cite{PR} \cite{PR2}, who computed the $K$-groups of
the twisted group
$C^*$-algebras
of uniform lattices in solvable groups. More precisely, they proved that
if $\Gamma$ is a uniform lattice in a solvable Lie group $G$, then
$$
K_*(C^*(\Gamma,\sigma)) \cong K^{* + \dim G}(\Gamma\backslash G,
\delta(B_\sigma))
$$
where $\sigma$ is any multiplier
on $\Gamma$, $K^{*}(\Gamma\backslash G,
\delta(B_\sigma)) $ denotes the twisted $K$-theory of a continuous trace
$C^*$-algebra
$B_\sigma$ with spectrum $\Gamma\backslash G$, while $\delta(B_\sigma)
\in H^3(\Gamma\backslash G, \mathbb Z)$ denotes the
Dixmier-Douady invariant of $B_\sigma$.  (Note that
the twisted $K$-theory was
studied in \cite{Ros}).
Packer and Raeburn proved a stabilization theorem and used the Thom isomorphism
theorem for the $K$-theory of $C^*$ algebras, due to Connes \cite{Co2}, to
prove their
results.

In Section 8 we extend the main theorem of \cite{PR}, \cite{PR2} to the
case when
$\Gamma$ is a lattice in a $K$-amenable Lie group $G$. More precisely, we
prove that
for such $G$ and $\Gamma$,
$$K_*(C^*(\Gamma,\sigma))\cong K_*(C^*_r(\Gamma,\sigma))$$
and
$$
K_*(C^*(\Gamma,\sigma)) \cong K^{* + dim(G/K)}(\Gamma\backslash G/K,
\delta(B_\sigma)),
$$
where $K$ is a maximal compact subgroup of $G$, $\sigma$ is any multiplier
on $\Gamma$,
$K^{*}(\Gamma\backslash G/K,
\delta(B_\sigma))$ is the twisted $K$-theory of a continuous trace $C^*$-algebra
$B_\sigma$ with spectrum $\Gamma\backslash G/K$, and $\delta(B_\sigma)\in
H^3(\Gamma\backslash G/K, \mathbb Z)$ is the
Dixmier-Douady invariant of $B_\sigma$.

Our method uses the $K$-amenability results of Kasparov \cite{Kas1} and the
Packer-Raeburn stabilization theorem \cite{PR}.
 In the case when $\Gamma =\Gamma_g$ is
the fundamental group of a Riemann surface $\Sigma_g$ of genus $g>0$, we
deduce that the Dixmier-Douady class $\delta(B_\sigma))$ is trivial.
Using this we demonstrate that  for any multiplier $\sigma$ on $\Gamma_g$
$$
K_0(C^*(\Gamma_g,\sigma)) \cong K^0(\Sigma_g) \cong \mathbb Z^2,\eqno
$$
and that
$$
K_1(C^*(\Gamma_g,\sigma)) \cong K^1(\Sigma_g) \cong \mathbb Z^{2g}.\eqno
$$
We end the discussion with an
interesting conjecture for compact 3-dimensional manifolds which are
Eilenberg-Maclane
spaces. These $K$-theoretic results have now been
generalized to $C^*$-dynamical systems in \cite{Ma}.

One of the most outstanding open problems about magnetic Schrodinger
operators or Hamiltonians
on Euclidean space is concerned with the nature of their spectrum,
and is called the {\em Ten Martini Problem} (TMP) (cf. \cite{Sh}). More
precisely, TMP asks whether
given a multiplier $\sigma$ on $\mathbb Z^2$, is there an associated Hamiltonian
({\it i.e.} a Hamiltonian which commutes with the $(\Gamma, \sigma)$ projective
action of $\Gamma$ on $L^2(\mathbb R^2)$) possessing a Cantor set
type spectrum, in the sense that the intersection of
the spectrum of the Hamiltonian with some compact interval in $\mathbb R$
is a Cantor set?
One can deduce from the range of the
trace on $K_0$ of the twisted group $C^*$-algebras that when the
multiplier takes its values in the  roots of unity in $U(1)$ (we say then that
it
is rational) that such a Hamiltonian cannot exist. However, in the Euclidean
case
and for Liouville numbers, the discrete analogue of the TMP has been been
settled
in the affirmative by Choi, Elliot and Yui  \cite{CEY} (cf. \cite{Sh} for a
historical
perspective). In Section 9 we are concerned also with the hyperbolic
analogue
of the TMP, which we call the {\em Ten Dry Martini Problem} (TDMP). We
prove that
the Kadison constant
 of the twisted group $C^*$-algebra $C^*_r(\Gamma_g,\sigma)$
is positive whenever the multiplier is rational, where $\Gamma_g$ is now
the fundamental
group of a genus $g$ Riemann surface. We then use the results of Br\"uning
and Sunada \cite{BrSu}
to deduce that when the multiplier is rational  the TDMP answered in the
 negative, and
we leave open the more difficult irrational case.
The calculation of the range of the trace
exploits a number of results including a
 twisted Kasparov map on $K$-theory.
 Finally, we
apply our results to  give a complete
classification up to
isomorphism of the twisted $C^*$-algebras $C^*_r(\Gamma_g,\sigma)$.

In  Sections 10 and 11 we will identify the
character of our Fredholm modules, the `Hall conductivity',
for both the continuous and discrete models. This character
$\tau_c(P,P,P) = \tr(P\, dP\, dP)$ is shown to arise
from Connes' `area cocycle' and we are able to
identify it with a topological invariant,
generalising the work of Xia \cite{Xia} in the case of
the quantum Hall effect on
Euclidean space. We use the pairing between $K$-theory and cyclic
cohomology \cite{Co}, a generalization of the Connes-Moscovici index
theorem \cite{CM} to projectively invariant elliptic operators
and the twisted analogue of the Kasparov map.
In fact we obtain a general index theorem
 which equates the (analytical) index arising from
the Fredholm modules to a topological index.
It specialises in the case of the
cyclic cocycle $\tau_c$ to give the
surprising fact that the hyperbolic `Hall conductivity'
$\tau_c(P,P,P)\in 2(g-1)\ints .$
This raises the obvious question of whether a real material
with a hyperbolic crystalline geometry could be manufactured
and the genus of the quotient Riemann surface measured experimentally.
To be specific the model we consider here can be understood most easily
in the imbedded
hyperboloid version of hyperbolic space. If we use the $\Gamma$-orbit
of a point in the hyperboloid in $\real^3$
 to represent the crystal lattice structure of a conducting
material
then our discrete model corresponds to applying a magnetic field
which is everywhere normal to the hyperboloid. (This captures the hyperbolic
geometry.) Then by regarding the lattice points as the vertices of a graph
whose edges are geodesics corresponding to the generators of $\Gamma$
our model Hamiltonian corresponds to allowing electrons to hop
between sites on the lattice along the edges of the graph.
Then our theorem predicts that the conductivity should depend on the genus
of the Riemann surface obtained by quotienting the hyperboloid by $\Gamma$.

In Section 6 we exhibit a cyclic cocycle which plays the role
of the Kubo formula for higher genus surfaces. It has an intrinsic
geometric description as a `symplectic area' cocycle on the
universal cover of the Jacobi variety of the Riemann surface.
The novel feature of the higher genus case
 (as opposed to genus one which is the
Euclidean case) is that the Kubo cocycle is cohomologous
(but not equal) to
the cyclic cocycle arising from the Fredholm module.
Given our $K$-theoretic interpretation of the latter this is sufficient
to give the anticipated result that the Hall conductivity, as defined
through the Kubo cocycle, is integral and depends on the genus.

We conclude by showing how our
formalism links with the non-commutative
Riemann surface theory described in  \cite{Klim+Les1}, \cite{Klim+Les2}.

\section{The geometry of the hyperbolic plane}

The upper half-plane can be mapped by the Cayley transform
$z = (\zeta-i)/(\zeta+i)$ to the unit disc $\disc$ equipped with the metric
$|dz|^2/(1-|z|^2)^2$ and symplectic form $dz\,d\overline{z}/2i(1-|z|^2)^2$, on
which $\su$ acts, and some calculations are more easily done in that setting.
In order to preserve flexibility we shall work more abstractly with a Lie
group $G$ acting transitively on a space $X \sim G/K$.
Although we shall ultimately be interested in the case of $G = \slr$ or $\su$,
and $K$ the  maximal compact subgroup which stabilises $\zeta = i$ or $z=0$,
those details will play little role in many of our calculations, though we
shall need to assume that $X$ has a $G$-invariant Riemannian metric and
symplectic form $\omega_\hyp$.
We shall denote by $\Gamma$ a discrete subgroup of $G$ which
acts freely on $X$ and hence  intersects $K$
trivially.

We shall assume that $\cl$ is a hermitian line
bundle over $X$, with a connection, $\nabla$, or equivalently, for each pair of
points $w$ and $z$ in $X$, we denote by $\tau(z,w)$ the parallel transport
operator
along the geodesic from $\cl_w$ to $\cl_z$.
In $\hyp$ with the line bundle trivialised and $\eta = \theta dx/y$ one can
calculate explicitly that
$$\tau(z,w) =   \exp\left(i\int_w^z\eta\right)
= [(z-\overline{w})/(w-\overline{z})]^\theta.$$
For general $\eta$ we have $\eta - \theta dx/y = d\phi$ and
$$\tau(z,w) = \exp(i\int_w^z\eta) =
[(z-\overline{w})/(w-\overline{z})]^\theta\exp(i(\phi(z)-\phi(w))).$$
Parallel transport round a geodesic triangle with vertices $z$, $w$, $v$,
gives rise to a holonomy factor:
$$\varpi(v,w,z) = \tau(v,z)^{-1}\tau(v,w)\tau(w,z),$$
and this is clearly the same for any other choice of $\eta$, so we may as well
work in the general case.

\begin{lemma}
The holonomy can be written as
$\varpi(v,w,z) = \exp\left(i\theta\int_\Delta \omega_\hyp\right)$,
where $\Delta$ denotes the geodesic triangle with vertices $z$, $w$ and $v$.
The holonomy is invariant under the action of $G$, that is
$\varpi(v,w,z) = \varpi(gv,gw,gz)$, and under cyclic permutations of its
arguments.
Transposition of any two vertices inverts $\varpi$.
For any four points $u$ ,$v$, $w$, $z$ in $X$ one  has
$$\varpi(u,v,w)\varpi(u,w,z) = \varpi(u,v,z)\varpi(v,w,z).$$
\end{lemma}
\begin{proof}
By definition, for a suitable trivialisation of $\cl$ one has
$$\varpi(v,w,z) = \exp\left(i\int_{\partial\Delta}\eta\right)$$
and the first part follows by applying Stokes' Theorem after noting that the
result is independent of the trivialisation.
The invariance under $G$ follows from the invariance of the symplectic form,
and the results of permutations follow from the properties of the integral,
as does the final identity.
\end{proof}

\section{The twisted algebra of kernels}

The geometrical data described in the last section enables us to easily
describe the first of the two C$^*$algebras which appear in the theory.
This twisted algebra of kernels, which  was introduced by Connes \cite{Co2}
is the C$^*$-algebra $\ck$ generated by compactly supported smooth
functions on
$X\times X$ with the multiplication
$$k_1*k_2(z,w) = \int_X k_1(z,v)k_2(v,w)\varpi(z,w,v)\,dv,$$
(where $dv$ denotes the $G$-invariant measure defined by the metric) and
$k^*(z,w)= \overline{k(w,z)}$.
There is an obvious trace on $\ck$ given by $\tau_\ck (k) = \int_X k(z,z)\,dz$
The algebra of twisted kernels is the extension of the C$^*$-algebra of the
principal groupoid $X\times X$ defined by the cocycle
$((v,w),(w,z)) \mapsto \varpi(v,w,z)$, \cite{Ren1}.

\begin{lemma}
The algebra $\ck$ has a representation $\pi$  on the space of $L^2$ sections
of $\cl$ defined by
$$(\pi(k)\psi)(z) = \int_X k(z,w)\tau(z,w)\psi(w)\,dw.$$
\end{lemma}
\begin{proof}
The parallel transport $\tau(z,w)$ ensures that the integral is in the
appropriate fibre, and the fact that it is a representation follows from
a calculation using the definition of the holonomy.
\end{proof}

Before describing the second algebra we need to link the geometrical data
more directly to the group $G$.
To do this we fix a basepoint $u\in X$ and introduce the function
$\phi$ from  $X\times G$ to line bundle automorphisms defined by
$$\phi(z,g) = \varpi(u, g^{-1}u,g^{-1}z)\tau(u,z)^{-1}\tau(u,g^{-1}z).$$
(The ratio of parallel transports defines an operator from the fibre
$\cl_{g^{-1}z}$ to $\cl_z$.)

\begin{lemma}
The function $\phi$ satisfies
\begin{eqnarray*}
\phi(z,x)\phi(x^{-1}z,y) &
= &  \varpi(u,y^{-1}u,y^{-1}x^{-1}u)\phi(z,xy)\\
\phi(z,x)\tau(x^{-1}z,x^{-1}w) & = & \tau(z,w)\phi(w,x).
\end{eqnarray*}
\end{lemma}
\begin{proof}
By definition we have
\begin{eqnarray*}
\phi(z,x)\phi(x^{-1}z,y) & = & \varpi(u,x^{-1}u,x^{-1}z)\varpi(u,
y^{-1}u,y^{-1}x^{-1}z)
\tau(u,z)^{-1}\\
& & \phantom{ijkijkkk}\tau(u,x^{-1}z)\tau(u,x^{-1}z)^{-1}\tau(u,y^{-1}x^{-1}z)\\
& = & \frac{\varpi(u,x^{-1}u,x^{-1}z)\varpi(u, y^{-1}u,y^{-1}x^{-1}z)}
{\varpi(u,y^{-1}x^{-1}u,y^{-1}x^{-1}z)}\phi(z,xy).
\end{eqnarray*}
Now by  Lemma 2.1
\begin{eqnarray*}
\lefteqn{\varpi(u,x^{-1}u,x^{-1}z)\varpi(u, y^{-1}u,y^{-1}x^{-1}z)}\\ &
\phantom{==============}= &
\varpi(y^{-1}u,y^{-1}x^{-1}u,y^{-1}x^{-1}z)\varpi(u,
y^{-1}u,y^{-1}x^{-1}z)\\
&\phantom{==============} = &
\varpi(u,y^{-1}x^{-1}u,y^{-1}x^{-1}z)\varpi(u, y^{-1}u,y^{-1}x^{-1}u),
\end{eqnarray*}
from which the first result follows.
For the second result we note (compressing the notation) that
\begin{eqnarray*}
\frac{\tau(x^{-1}z,x^{-1}w)\phi(z,x)}{\phi(w,x)} & = &
 \frac{\varpi(u, x^{-1}u,x^{-1}z)}{\varpi(u, x^{-1}u,x^{-1}w)}
\frac{\tau(x^{-1}z,x^{-1}w)\tau(u,x^{-1}z)}{\tau(u,z)}
\frac{\tau(u,w)}{\tau(u,x^{-1}w)}\\
 & = & \frac{\varpi(u, x^{-1}u,x^{-1}z)\varpi(u,x^{-1}z,x^{-1}w)}
{\varpi(u, x^{-1}u,x^{-1}w)}
\frac{\tau(u,x^{-1}w)}{\tau(u,z)}
\frac{\tau(u,w)}{\tau(u,x^{-1}w)}\\
 & =  & \varpi(x^{-1}u,x^{-1}z,x^{-1}w)\frac{\tau(u,w)}{\tau(u,z)}\\
& = & \varpi(u,z,w)\frac{\tau(u,w)}{\tau(u,z)}\;
= \;\tau(z,w).
\end{eqnarray*}
\end{proof}

The most important  aspect of the first result is that
$$\sigma(x,y) = \phi(z,xy)/\phi(z,x)\phi(x^{-1}z,y)
= \varpi(u,y^{-1}u,y^{-1}x^{-1}u)^{-1}
= \varpi(u,xu,xyu)$$
is independent of $z$.
(We note also that $\sigma(g,1) = \sigma(1,g) = \sigma(g,g^{-1}) = 1$.
Although these normalisations do not seriously affect matters they can
sometimes be used to simplify formulae.)

\begin{lemma}
The function $\sigma: G\times G \to \torus$ satisfies the cocycle identity,
$$\sigma(x,y)\sigma(xy,g) = \sigma(x,yg)\sigma(y,g).$$
\end{lemma}
\begin{proof}
This is a simple calculation along the lines of those above.
\end{proof}

This result means that $\sigma$ defines a projective multiplier
or group 2-cocycle for $G$,
moreover, it is clearly continuous
and identically 1 when restricted  to $G\times K$ and to $K\times G$.

\begin{lemma}
The group G has a natural unitary $\sigma$-representation $U$ on the $L^2$
sections of $\cl$ defined by
$$U(g)\psi(z) = \phi(z,g)\psi(g^{-1}z).$$
\end{lemma}
\begin{proof}
This follows immediately from  lemma  3.4.
\end{proof}

This projective representation induces an action of $G$ as automorphisms of
$\ck$.

\begin{lemma}
For any $g\in G$ and $k \in \ck$ we have
$U(g)\pi(k)U(g)^{-1} = \pi(g.k)$, where
$$g.k(z,w) = k(g^{-1}z,g^{-1}w).$$
\end{lemma}
\begin{proof}
By direct calculation and use of Lemma 3.2
$$(U(g)\pi(k)U(g)^{-1}\psi)(z)
= \int \phi(z,g)k(g^{-1}z,g^{-1}w)\tau(g^{-1}z,g^{-1}w)
\phi(w,g)^{-1}\psi(g^{-1}w)\,dw$$
$$= \int k(g^{-1}z,g^{-1}w)\tau(z,w)\psi(g^{-1}w)\,dw,$$
from which the result follows.
\end{proof}

The second part of Lemma 3.2 can now be interpreted as saying that the parallel
transport $\tau$ behaves covariantly under $U(g)$, that is conjugation by
$U(g)$ sends $\tau(z,w)$ to $\tau(g^{-1}z,g^{-1}w)$.
Taking $w = \exp(-tX)z$ and considering the limit as $t \to 0$ we obtain the
following result:

\begin{cor}
The $\sigma$-representation $U$ and connection $\nabla$ are related by
$U(g)\nabla U(g)^{-1} = g.\nabla$,
where $g.\nabla$ denotes the natural action of $G$ on forms.
\end{cor}

\section{Various $C^*$-algebras}

\subsection{The imprimitivity algebra}
The $\sigma$-representation $U$ defined in the previous section is
clearly equivalent to one induced from a $\sigma$-representation, $L$, of the
isotropy subgroup $K$.
Such representations are characterised by the fact that they also admit an
action of the imprimitivity algebra.
In general this can be defined as one of Green's twisted crossed product
C$^*$-algebras \cite{Green} , but in the case of a continuous multiplier
$\sigma$ there is a
simpler direct construction.
The imprimitivity algebra, $\ca = \ca(G,K,\sigma)$, on $G/K$ is a completion
of the algebra $\ca_0 = C_c(G/K,G)$ with multiplication
$$(\alpha*\beta)(s,g) =
\int_G \alpha(s,x)\beta(x^{-1}s,x^{-1}g)\sigma(x,x^{-1}g)^{-1}\,dx$$
and involution
$$\alpha^*(s,g) = \sigma(g,g^{-1})\overline{\alpha(g^{-1}s,g^{-1})}.$$
(With the conventions of the last section $\sigma(g,g^{-1})=1$ and could be
omitted.)
These formulae use the unimodularity of $G$ and the existence of a
$G$-invariant measure on $G/K$, otherwise some Radon-Nikodym derivatives
would be needed.

The algebra has a trace
$$\tr_\ca(\alpha) = \int_{G/K} \alpha(s,1)\,ds.$$
More details may be found in \cite{Green}, where it is also shown that
$\ca(G,K,\sigma)$ is Morita equivalent to
$C^*(K,\sigma)\otimes{\mathcal K}(L^2(G/K))$,
where $C^*(K,\sigma)$ denotes the twisted group $C^*$-algebra.
Like $\ck$ the imprimitivity algebra is a groupoid algebra (being an extension
by $\sigma$ of the algebra of the transformation groupoid for $G$ acting on
$X$) and most of this paper could be understood in the context of groupoids,
\cite{Ren1}.

The algebra $\ck$ could have been derived from $\ca$ as a quotient, as we
shall now show.
Let $L$ be a $\sigma$-representation of $K$ on a Hilbert space $\chh_L$.
(Since the multiplier of the last section is 1 whenever either of its
arguments is in $K$, one could take $L = 1$, but the argument works more
generally.)
For each $L$ the imprimitivity  algebra has a natural $*$-representation on
the induced representation space of $\chh_L$-valued functions on $G$ satisfying
the equivariance condition
$$\psi(gk) = \sigma(g,k)^{-1}L(k)^*\psi(g)$$
for all $g\in G$ and $k\in K$. This representation is given by
$$(\alpha.\psi)(z) =
\int_G \alpha(zK,x)\sigma(x,x^{-1}z)^{-1}\psi(x^{-1}z)dx.$$
(It may be checked that $\alpha.\psi$ satisfies the same equivariance
condition as $\psi$.)
The group $G$ has an induced $\sigma$-representation on this function space
given by
$$U(g)\psi(z) = \sigma(g,g^{-1}z)^{-1}\psi(g^{-1}z).$$
The imprimitivity algebra incorporates both this action and the multiplication
operators, and so permits the description of quantum mechanical momentum and
position operators on $G/K$.
The group action allows for the free Hamiltonian $-\frac 18{\bf J}.{\bf J}$,
whilst the functions on $X = G/K$ make it possible to
add  an extra potential, $V$.

For an appropriate choice of $L$, $U$ is equivalent to the representation
in the last
section.
Indeed we may identify the equivariant functions on $G$ with sections of
the line bundle $\cl$ and then we have, in the previous notation,
$$(\alpha.\psi)(z) =
\int_G \alpha(zK,x)\phi(zK,x)\psi(x^{-1}z)\,dx.$$

In this form we may easily see the connection to the algebra $\ck$ of twisted
kernels.

\begin{lemma}
There is a $*$-homomorphism $\alpha \mapsto T_\alpha$, defined by
$$T_\alpha(z,w) =
\int_G \frac{\alpha(z,x)\phi(z,x)}{\tau(z,w)}\delta(w,x^{-1}z)\,dx$$
from $\ca$ onto $\ck$
(where $\delta$ just restricts the integration to those $x$ satisfying $w =
x^{-1}z$).
\end{lemma}
\begin{proof}
We calculate that
\begin{eqnarray*}
\lefteqn{(T_\alpha*T_\beta)(z,w)} \\
& = & \int T_\alpha(z,v)T_\beta(v,w)\frac{\varpi(z,v,w)}{\tau(z,v)\tau(v,w)}
\,dv\\
& = & \int \alpha(z,x)\beta(v,x^{-1}y)\frac{\phi(z,x)\phi(v,x^{-1}y)}{\tau(z,w)}
\delta(v,x^{-1}z)\delta(w,y^{-1}xz)\,dx\,dy\,dv\\
& = & \int \alpha(z,x)\beta(x^{-1}z,x^{-1}y)
\frac{\phi(z,x)\phi(x^{-1}z,x^{-1}y)}{\tau(z,w)}
\delta(v,x^{-1}z)\delta(w,y^{-1}z)\,dx\,dy\,dv\\
& = & \int \alpha(z,x)\beta(x^{-1}z,x^{-1}y)
\frac{\phi(z,y)}{\sigma(x,x^{-1}y)\tau(z,w)}
\delta(v,x^{-1}z)\delta(w,y^{-1}z)\,dx\,dy\,dv\\
& = &  \int (\alpha*\beta)(z,y)\frac{\phi(z,y)}{\tau(z,w)}\delta(w,y^{-1}z)\,dy
\; = \; T_{\alpha*\beta}(z,w),
\end{eqnarray*}
which proves the homomorphism property.
It follows similarly that it is a  $*$-homomorphism.

To see that it is surjective we note that if the Haar measure on $K$ is
normalised then
$$\alpha(z,x) = T_\alpha(z,x^{-1}z)\frac{\tau(z,x^{-1}z)}{\phi(z,x)}$$
provides an inverse.
(The key is to note that if $g(z)u = z$ and $g(w)u = w$ then $w=x^{-1}z$
forces $x$ to have the form $g(z)kg(w)^{-1}$ for some $k \in K$.)
\end{proof}

This may be interpreted as saying that $\ck$ is a quotient of $\ca$, and this
means that the behaviour of the two algebras is very similar.
For this reason we shall often merely show the constructions in the case of one
and leave it as an exercise to fill in the details for the other.
However, it will be useful to note that $\ck$ has a natural trace
$$\tr_\ck(T) = \int_X T(z,z)\,dz.$$

\subsection{The regular representation}
The trace gives rise to an inner product
$$\ip{\alpha}{\beta} = \tr_\ca(\alpha^**\beta)
= \int_{G/K\times G} \overline{\alpha(s,x)}\beta(s,x)\,ds\,dx$$
on $\ca_0$, and completion with respect to this gives rise to a Hilbert space
$\chh_\ca$ on which $\ca_0$ is represented by left multiplication.
Since the action is continuous this extends to give the left regular
representation of $\ca$.
This is given by the same formula as the algebra multiplication. That is, for
$\Psi\in \chh_\ca$, one has
$$(\alpha .\psi)(s,g) =
\int_G \alpha(s,x)\psi(x^{-1}s,x^{-1}g)\sigma(x,x^{-1}g)^{-1}\,dx.$$
There is also a $\sigma$-representation of $G$ on $\chh_\ca$, obtained by
sending $g \in G$ to the function $(s,x) \mapsto \delta_g(x)$.

For any unitary character $\chi$ of $K$, one may define a generalised
regular representation on the space $H_\ca^\chi$ of functions
$\psi \in C_c(G\times G)$ which satisfy the condition
$$\psi(zk,g) = \chi(k)\psi(z,g).$$
Since $|\psi(zk,g)|^2$ is independent of $k$, we may use the same inner
product as before, and it is easy to check that the action
$$(\alpha .\psi)(z,g) =
\int_G \alpha(zK,x)\psi(x^{-1}z,x^{-1}g)\sigma(x,x^{-1}g)^{-1}\,dx$$
respects the equivariance condition.
When $\chi = 1$ we obtain the regular representation.

We could similarly define the regular representation of $\ck$ and also a
generalisation defined by
$$(T.\Psi)(z,w) = \int_X T(zK,v)\Psi(v,w)\,dv$$
on the space of kernels $\Psi$ on $G\times G$ which satisfy
$\Psi(zk_1,wk_2) = \chi(k_2)\Psi(z,w)$.

\subsection{The $\Gamma$-invariant imprimitivity algebra}
Let $\Gamma$ be another subgroup of $G$, and let $\ca^{\Gamma}$ denote the
part of the imprimitivity algebra which commutes with the induced
representation of $\Gamma$.
This is readily seen to consist of those functions $\alpha \in \ca(G,K,\sigma)$
which satisfy the condition $\gamma.\alpha = \alpha$, where
$$(\gamma.\alpha)(s,g) =
\sigma(g,\gamma)\sigma(\gamma,\gamma^{-1}g\gamma)^{-1}
\alpha(\gamma^{-1}s,\gamma^{-1}g\gamma)$$
for all $\gamma \in \Gamma$.
To see this we note that
$$(\delta_\gamma*\alpha)(s,g) =
\sigma(\gamma,\gamma^{-1}g)\alpha(\gamma^{-1}s,\gamma^{-1}g),$$
which is identical to
$$(\gamma.\alpha*\delta_\gamma)(s,g)
= \sigma(g\gamma^{-1},\gamma)^{-1}\gamma.\alpha(s,g\gamma^{-1})
= \sigma(\gamma,\gamma^{-1}g)^{-1}\alpha(\gamma^{-1}s,\gamma^{-1}g).$$

When $g=1$ the condition $\gamma.\alpha(s,g) = \alpha(s,g)$ reduces to
$\alpha(\gamma^{-1}s,1) = \alpha(s,1)$ so that one obtains a trace
$$\tr_{\ca^\Gamma}(\alpha) = \int_{\Gamma\backslash G/K} \alpha(s,1)\,ds.$$
(More generally, the $\Gamma$-invariant functions are determined by
their values at a single point $s$ of each $\Gamma$ orbit on $G/K$.)
Using $\tr_{\ca^\Gamma}$ one may define a regular representation of $\cag$.

In the case of the twisted kernel algebra one may likewise pick out a
$\Gamma$-invariant subalgebra $\ck^\Gamma$, which commutes with $U(\gamma)$
for all $\gamma\in \Gamma$.
Using Lemma 3.2 this condition reduces simply to the requirement that the
kernel satisfies $k(\gamma^{-1}z,\gamma^{-1}w) = k(z,w)$.
The natural trace $\tr_{\ck^\Gamma}$ for this algebra is given by the same
formula as before
except that the integration is over a fundamental domain $X_\Gamma$
 rather than $X$:
$$
\tau_{\ck^\Gamma} (T)=\int_{X_\Gamma} T(z,z)dz.
$$

\subsection{Morita equivalence}
Later we shall need some $K$-theory, and so it will be useful to show
that the algebra $\ck^\Gamma$ is Morita equivalent to another more tractable
algebra.
We shall do this by using the groupoid equivalence arguments of \cite{M+R+W}, or
rather the twisted version, \cite{Ren2}, \cite{Ren3}.
We have already noted that $\ck$ is an extension of the groupoid $X\times X$ by
a cocycle defined by $\varpi$, and $\Gamma$ invariance of $\varpi$ means that
$\ck^\Gamma$ is likewise the extension of $X\times_\Gamma X$ by $\varpi$, where
$X\times_\Gamma X$ denotes the groupoid obtained by factoring out the diagonal
action of $\Gamma$.
More precisely, the groupoid elements are $\Gamma$ orbits
$(x,y)_\Gamma = \{(\gamma x,\gamma y): \gamma \in \Gamma\}$,
and $(x_1,y_1)_\Gamma$ and $(x_2,y_2)_\Gamma$ are composable if and
only if $y_1 = \gamma x_2$ for some $\gamma \in \Gamma$, and then the
composition is $(x_1,\gamma y_2)_\Gamma$.

\begin{thm}
The algebra $\ck^\Gamma$ is Morita equivalent to the twisted group algebra
$C^*(\Gamma,\overline{\sigma})$.
\end{thm}
\begin{proof}
This result will follow immediately from \cite{Ren2} Corollaire 5.4 (cf
\cite{M+R+W}
Theorem 2.8) once we have established the groupoid equivalence in the
following Lemma.
\end{proof}
\begin{lemma}
The line bundle $\cl$ over $X$ provides an equivalence (in the sense of
\cite{Ren2}
Definition 5.3) between the groupoid extensions $(X\times_\Gamma X)^\varpi$ of
$X\times_\Gamma X$ defined
by $\varpi$ and $\Gamma^\sigma$ of $\Gamma$ defined by $\overline{\sigma}$.
\end{lemma}
\begin{proof}
Both extensions are by $\torus$.
We write the elements of $(X\times_\Gamma X)^\varpi$ as triples
$(x,y,t) \in X\times X\times \torus$ with the first two
elements representing a diagonal $\Gamma$ orbit.
Elements are composable if their first two components are composable, and,
when $y_1=\gamma x_2$,
$$(x_1,y_1,t_1)(x_2,y_2,t_2) =
(x_1,\gamma y_2, t_1t_2\varpi(x_1,y_1,\gamma y_2)).$$

The line bundle can be trivialised and written as $X\times \comp$.
We let $(X\times_\Gamma X)^\varpi$ act on the left of the line bundle by
defining $(x,y,t)$ to act on $(z,u)$ if $z = \gamma y$ for some
$\gamma \in \Gamma$, and then the result of the action is
$(\gamma x,\tau(\gamma x,z)tu)$.
(One may check that this gives an action using the relationship between
parallel transport and holonomy and the $\Gamma$-invariance of $\varpi$.)

The twisted groupoid $\Gamma^\sigma$ has underlying set $\Gamma\times \torus$,
all elements are composable, and multiplication is given by
$$(\beta,s)(\gamma,t) = (\beta\gamma,\sigma(\beta,\gamma)st).$$
It acts on the right of $\cl$ by
$$(z,u).(\gamma,t) = (\gamma^{-1}z,\phi(z,\gamma)^{-1}tu).$$
(The fact that this defines an action follows from the definition of $\sigma$
in terms of $\phi$.)
We may now check that these actions commute, since, if $z=\beta y$,
$$[(x,y,t)(z,u)](\gamma,s) = (\beta x,\tau(\beta x,z)tu)(\gamma,s)
= (\gamma^{-1}\beta x,\phi(\beta x, \gamma)^{-1}\tau(\beta x,z)tus),$$
whilst
$$(x,y,t)[(z,u)(\gamma,s)] = (x,y,t)(\gamma^{-1}z,\phi(z,\gamma)^{-1}us)
= (\gamma^{-1}\beta x,
\tau(\gamma^{-1}\beta x,\gamma^{-1}z)\phi(z,\gamma)^{-1}tus),$$
and the equality of these two follows from the second part of Lemma 3.2.
\end{proof}

\begin{rems*}
We will observe in Section 9 that  the algebra
$C^*(\Gamma,\overline{\sigma})$
is isomorphic to
$C^*(\Gamma,{\sigma})$
which in turn is known to be isomorphic to the imprimitivity
algebra for $\sigma$-inducing from $\Gamma$ to
 $C^*(G,{\sigma})$. This latter algebra, denoted  $C^*(G/\Gamma,G)$ is
the completion of   $C_c(G/\Gamma,G)$ where the
latter has a multiplication analogous to that described above
for  $C_c(G/K,G)$ (simply replace $K$ by $\Gamma$ in the earlier
discussion). Thus a corollary of our results in this
subsection is that $\ck^\Gamma$ is Morita equivalent
to  $C^*(G/\Gamma,G)$. Furthermore our discussion below of a Fredholm module
for   $\ck^\Gamma$ may be modified so as to produce a Fredholm module
for  $C^*(G/\Gamma,G)$ whose character, for $G=\so(n,1)$,
 is also given by the area
cocycle. We omit the details here as they would take us too far afield
(see however the analogous discussion in the discrete case in \cite{Co2}).
\end{rems*}

\subsection{The Hamiltonian}

We have asserted informally that the Hamiltonian can be accommodated within the
algebras $\ca^\Gamma$ and $\ck^\Gamma$ and we shall now provide the proof.
We work with the smaller algebra $\ck^\Gamma$, the results for $\ca^\Gamma$
following similarly.

\begin{lemma}
The Hamiltonian $H = -\frac 12\nabla^*.\nabla$ commutes with the projective
representation $U$.
\end{lemma}
\begin{proof}
We recall from Corollary 3.6 that $U(g)\nabla U(g)^{-1} = g.\nabla$, so that
$$U(g)HU(g)^{-1} = -\frac 12(g.\nabla)^*(g.\nabla),$$
and, since the Riemannian structure is invariant under the action of $G$, this
is just $H$.
\end{proof}

We could also obtain the same result by writing
$H = \frac 18{\bf J}.{\bf J}+\frac 14B^2$ (for some constant $B$)
 and using the fact that
the Casimir operator commutes
with the representation, provided that we check that $J_k$ are the
representatives of the Lie algebra generators in the representation $U$.
Using the invariance of $H$ it is now not difficult to see the following
result:

\begin{lemma}
The Hamiltonian $H$ is affiliated to the von Neumann algebra generated by the
representation $\pi$ of $\ck$.
\end{lemma}
\bigskip
Since $H$ is $G$-invariant and so {\it a fortiori} also $\Gamma$-invariant, it
is sufficient to look at $\ck^\Gamma$.

\begin{cor}
The Hamiltonian $H$ is affiliated to the von Neumann algebra generated by the
representation $\pi$ of $\ck^\Gamma$.
\end{cor}
\bigskip
We next observe that Br\"uning and Sunada have proved an estimate on
the Schwartz kernel of $\exp(-tH)$ for $t > 0$, which implies that it is $L^1$
in each variable separately.
Since this kernel is $\Gamma$-invariant (by Lemma 7.1) it follows (in exactly
the same fashion as Lemma 4 of \cite{BrSu}) that this estimate implies that
$\exp(-tH)$ is actually in the algebra $\ck^\Gamma$.

\begin{lemma}
The operator $e^{-tH}$ is an element of $\ck^\Gamma$.
\end{lemma}
\begin{cor}
The spectral projections of $H$ corresponding to gaps in the spectrum lie in
$\ck^\Gamma$.
\end{cor}
\begin{proof}
If $\mu$ lies in a gap of the spectrum of $H$ then let $f$
be a continuous approximate step function which is identically one
on the part of the spectrum of $H$ contained in $[0,\mu]$ and zero
on the part contained in $[\mu,\infty)$.
Define $g(x)=f\circ\ln(1/x)$ for $x\in [0,1]$. Then $g$ is
a bounded continuous function which, when applied to $e^{-tH}$
gives the spectral projection corresponding to the interval $[0,\mu]$.
\end{proof}
\bigskip
Finally we consider the interacting Hamiltonian $H+V$, where $V$ is a
$\Gamma$-invariant function on $X$.
Notice that if $\psi$ is a continuous function of compact support on
$(X\times X)/\Gamma$ then
$V\psi$ is also such a function and hence defines an element
in the groupoid algebra.
Now, by Lemma 11, the resolvent of $H$ lies in $\ck^\Gamma$ and by
writing
$$(z-H-V)^{-1} = (1-(z-H)^{-1}V)^{-1}(z-H)^{-1},$$
and expanding $(1-(z-H)^{-1}V)^{-1}$ in a power series, we see that the
resolvent of $H+V$ is in the algebra $\ck^\Gamma$. This entails by
a simple modification of Corollary 3
that the spectral projections of $H+V$
corresponding to a gap in the spectrum also lie in $\ck^\Gamma$.

\section{ The discrete model}

In this section we formulate a version of the integer
quantum Hall effect on a graph   in hyperbolic space.
 The discussion uses a construction due to Sunada \cite{Sun}
 together
with a modification of Connes' work on Fredholm modules for the
group $C^*$-algebra of a discrete subgroup of a Lie group
\cite{Co2}.

The graph  is obtained by taking our group $\Gamma$,
the fundamental group of a Riemann surface, which acts freely on
hyperbolic space, fixing a base point $u$ and taking the orbit through
$u$ under the $\Gamma$ action. This gives the vertices of the graph. The
edges of the graph are geodesics constructed as follows.
Each element of the group
may be written as a word of minimal length in the
$2g$ generators and their inverses. Each generator
and its inverse determines a unique geodesic emanating from a vertex $x$
and these form the edges  of the
graph. Thus each word $x$ in the generators determines a
piecewise geodesic path from $u$ to $x$.

Sunada constructs a Hamiltonian on $\ell^2(\Gamma.u)$
which is a generalised Harper operator.
This construction, specialised to our case,
has the following form.
First we note that we may trivialise the restriction of the line bundle
to the vertices  and so
without loss of generality
 the appropriate Hilbert space becomes $\ell^2(\Gamma.u)$.
While the construction works for
 any connection 1-form $A$ on a line bundle $L$ over hyperbolic
space we make the formulae explicit by restricting
to the case where $A$ is the one form $\eta$.
For each directed edge $e$ of the graph joining $o(e)$ to $t(e)$
we define a function
$\tau(e)=\exp(i\int_e \eta)$.
Then $\tau(e)$ satisfies
$$\tau(\gamma.e)=\tau(e)(\frac{c\overline{o(e)}+d}{c{o(e)} +d})^\theta
(\frac{c{t(e)}+d}{c\overline{t(e)}+d})^\theta$$
where $\gamma = \bigl(\begin{smallmatrix}
a & b \\ c & d
\end{smallmatrix}\bigr)$.
We introduce the notation
$$s_\gamma(w)=(\frac {c\bar{w}+d}{cw +d})^\theta.$$

\begin{rems} Note that in our earlier notation these definitions
amount to
$$\tau(e)\equiv\tau(t(e),o(e)),$$
while the function from $\Gamma\times H$ to $U(1)$
given by $(\gamma, w)\mapsto s_\gamma(w)$
is a projective 1-cocycle for the $\Gamma$ action on $H$ which
is cohomologous to  the function mapping
$(\gamma, w)\mapsto \phi(\gamma.w,\gamma)$.
\end{rems}

Consequently there is a projective action of
$\Gamma$ on $\ell^2(\Gamma.u)$
given on $f\in\ell^2(\Gamma.u)$ by
$$\tilde \rho_\gamma f(w)
 = s_\gamma(\gamma^{-1}w)f(\gamma^{-1}w).$$
We have
$$\tilde\rho_{\gamma_1\gamma_2}=\frac{s_{\gamma_1\gamma_2}(\gamma_2^{-1}
\gamma_1^{-1}x)}{s_{\gamma_1}(\gamma_1^{-1}x)s_{\gamma_2}(\gamma_2^{-1}
\gamma_1^{-1}x)}\tilde\rho_{\gamma_1}\tilde\rho_{\gamma_2},$$ and Sunada
shows that the function multiplying $\tilde\rho_{\gamma_1}
\tilde\rho_{\gamma_2}$ is independent of $x$. This is exactly the relation
found at the end of lemma 3: the correspondence is given
by $\phi(z,\ast)\leftrightarrow s_{\ast}(\ast^{-1}z)$.
To obtain an explicit expression for this function, let
$\gamma_3=\gamma_1\gamma_2$, and write
$(\gamma_i)=\bigl(\begin{smallmatrix}
a_i & b_i \\ c_i & d_i
\end{smallmatrix}\bigr)$ for $i=1,2,3$.
Then a direct calculation shows that $$\tilde\rho_{\gamma_1\gamma_2}=\sigma
(\gamma_1,\gamma_2)\tilde\rho_{\gamma_1}\tilde\rho_{\gamma_2},$$ where
$\sigma$ is given by
$$\sigma(\gamma_1,\gamma_2)=\exp (2i\theta (\arg{(c_3\overline{w}+d_3)}-
\arg{(c_1(\overline{\gamma_2w})+d_1)}-\arg{(c_2\bar w+d_2)})).$$
It is not difficult to see that the right hand side is independent of the
choice of $w$.

Following \cite{Sun}, we
define our discrete hamiltonian, for $f\in\ell^2(\Gamma.u)$, by
$$h_\tau f(w) = \sum_{\gf{e}{o(e)=w}}\!\!\tau(e)f(t(e)).$$
Then $h_\tau$ is a generalised difference operator  (Sunada shows
that the Harper operator arises in a similar fashion).
One can verify by direct calculation that $h_\tau$
commutes with the projective action of $\Gamma$. This will, however,
become readily apparent when we
 transfer this construction to
$\ell^2(\Gamma)$. Define $\iota:\ell^2(\Gamma.u)\rightarrow\ell^2(\Gamma)$
by:
$$\iota(f)(\gamma)=\tilde\rho(\gamma) f(u).$$
Observe that
$$\iota(\tilde\rho(\mu)f)(\gamma)=\sigma(\gamma,\mu)\iota(f)(\gamma\mu).$$
Thus $\iota$ intertwines $\tilde\rho$ with the $\bar\sigma$-representation:
$$\rho(\mu)\iota(f)(\gamma)=\sigma(\gamma,\mu)\iota(f)(\gamma\mu).$$
Henceforth we use $\iota$ to identify $\ell^2(\Gamma.u)$ with $\ell^2(\Gamma)$.

\begin{prop}\cite{Sun}
The operator $h_\tau$ on  $\ell^2(\Gamma.u)$
maps to the operator $H_\tau$ on $\ell^2(\Gamma)$ under $\iota$,
where
$$H_\tau f(\gamma)= \sum_{\mu\in\Gamma}\sigma(\mu,\mu^{-1}\gamma)^{-1}a(\mu)
f(\mu^{-1}\gamma)$$
and $a$ is the function on $\Gamma$ given by
$$a(\gamma)=
s_\gamma(u)\!\!\sum_{\gf{e}{\gf{o(e)=u}{t(e)=\gamma.u}}}\!\!\tau(e).$$
\end{prop}

\begin{cor}
The bounded self-adjoint operator
$H_\tau$ is in  the algebraic twisted group algebra $\comp(\Gamma,\sigma)$
(the elements of
finite support in $C^*(\Gamma,\sigma)$)
 as the function $a$ has finite support. Moreover $H_\tau$ acts on the left
as an element of this algebra and so
 commutes with the $\bar\sigma$-representation $\rho$, as the
 latter acts on the right.
 \end{cor}

\def\tr{{\rm tr}}

\section{ A Kubo formula}


\subsection{Conductivity cocycles}
In this subsection we present  an argument which
derives analogues of the Kubo formula for the
hyperbolic `Hall conductivity'. In subsection 7.2
we show how to construct, from the results of this subsection,
 a unique cocycle which may be
compared with the character of the Fredholm module of Section 7.

Our reasoning here is that
 the Hall conductivity in the Euclidean situation
is measured experimentally by determining the equilibrium
 ratio of the current in the direction of the applied electric field to the
 Hall voltage, which is the potential difference in the orthogonal direction.
 To calculate this mathematically we instead determine
 the component of the induced current that is orthogonal to the applied
 potential. The conductivity can then be obtained by dividing this
 quantity by the magnitude of the applied field.
In the hyperbolic case it would seem
 at first sight that there are no preferred directions.
However interpreting the generators of the fundamental group
as geodesics on hyperbolic space gives a family of preferred directions
emanating from the base point. For each pair of
directions it is therefore natural to imitate the procedure of the
Euclidean case
and mathematically this is done as follows.

The Hamiltonian $H$ in a magnetic field depends on the magnetic vector
potential ${\bf A}$ and the functional derivative $\delta_kH$
of $H$ with respect to the components of ${\mathbf A}$, denoted
$A_k$, gives
 the current density $J_k$.
(For simplicity we take variations within a one-parameter family.)
The expected value of the current in a state described by a projection
operator $P$ is therefore $\tr(P\delta_kH)$
(cf \cite{Av+S+Y} equation (3.2)).
The following claim is not proved by a rigorous argument: one needs to
check various analytical details as in \cite{Xia}.
We have refrained from doing so here as this
would take us too far afield from the main point,
namely, obtaining a sensible hyperbolic Kubo formula
which may be compared with the character of the Fredholm module
constructed in the next section.
For the moment $\tr$ will denote a generic trace. We will become specific
after extracting a rigorous definition of the Kubo formula.

\noindent{\bf Claim}:
$$\tr(P\delta_kH) = i\tr(P[\partial_tP,\delta_kP]).$$

\noindent{\bf Plausability argument}:
By using the invariance of the trace under the adjoint action of operators
and the equation of motion we see that
\begin{align*}
\tr(P[\partial_tP,\delta_kP]) &= -\tr([P,\delta_kP]\partial_tP)\\
&= -i\tr([P,\delta_kP][P,H])\\
&= i\tr([P,[P,\delta_kP]]H).
\end{align*}
Now $\delta_kP = \delta_k(P^2) = P(\delta_kP) + (\delta_kP)P$, whence
$P(\delta_kP) P = 0$ and we have
\begin{align*}
[P,[P,\delta_kP]] &= P(P(\delta_kP) - (\delta_kP)P)
-(P(\delta_kP) - (\delta_kP)P)P\\
&= P(\delta_kP) +(\delta_kP)P = \delta_kP.
\end{align*}
Consequently we may write
$$\tr(P[\partial_tP,\delta_kP]) = i\tr((\delta_kP)H)
= i\tr(\delta_k(PH)) - i\tr(P(\delta_kH)),$$
and, assuming that the trace is invariant under variation of $A_k$, the
first term vanishes leaving the result asserted.

\bigskip
If the only $t$-dependence in $H$ and $P$ is due to the variation of $A_j$, a
component distinct from $A_k$, then
$\partial_t = \partial A_j/\partial t\times\delta_j$.
Working  in the Landau gauge so that the electrostatic
potential vanishes, the electric field is given by
${\bf E} = -\partial {\bf A}/\partial t$, and so $\partial_t = -E_j\delta_j$.
Combining this with the previous argument we arrive at the following result:

\begin{cor}
The conductivity for currents in the $k$ direction induced by electric fields
in the $j$ direction is given by
$-i\tr(P[\delta_jP,\delta_kP])$.
\end{cor}
\begin{proof}
The expectation of the current $J_k$ is given by
$$\tr(P\delta_kH) = i\tr(P[\partial_tP,\delta_kP])
= -iE_j\tr(P[\delta_jP,\delta_kP]),$$
from which the result follows immediately.
\end{proof}\subsection{The derivations on a Riemann surface}
On a Riemann surface it is natural to investigate changes in the potential
corresponding to adding multiples of the real and imaginary parts of
holomorphic 1-forms.
(For the genus 1 torus with imaginary period this amounts to choosing
forms whose integral round one sort of cycle vanishes but the
 integral round the
other cycle is non-trivial.
This corresponds to putting a non-trivial voltage across one cycle and
measuring a current round the other.)

Let $\Sigma_g = \hyp/\Gamma$ be the Riemann surface determined by
quotienting by $\Gamma$.
We follow the usual conventions (see for example \cite{GH}) in
 fixing representative
homology generators corresponding to cycles $A_j,B_j, j=1,2,\ldots,g$
with each pair $A_j, B_j$ intersecting in a common base point and
all other intersection numbers being zero. We let
$a_j,j=1,2,\ldots,2g$
be harmonic 1-forms dual to
this homology basis (this means that
$a_j, j=1,\ldots,g$ are dual to $A_j, j=1,\ldots,g$
and $a_{j+g}, j=1,2,\ldots,g$ are dual to
$B_j,j=1,2,\ldots,g$).
\begin{defn*}
Let
$$\delta_j\tau(z,w) = i\int_w^z \alpha_j\,\tau(z,w)$$
and
$$\delta_j\varpi(u,g^{-1}u,g^{-1}z) =
i\int_{\partial\Delta} a_j\,\varpi(u,g^{-1}u,g^{-1}z)$$
where $\Delta$ is a triangle with vertices at the three arguments of $\varpi$.
\end{defn*}
One then calculates that
\begin{eqnarray*}
\delta_j\phi(z,g) &  = &
i\left(\int_{\partial\Delta} a_j -\int_u^z a_j + \int_u^{g^{-1}z}a_j
\right)\phi(z,g)\\
& = & i\left(\int_u^{g^{-1}u} a_j + \int_{g^{-1}u}^{g^{-1}z}a_j
- \int_u^z a_j\right)\phi(z,g).
\end{eqnarray*}
This can also be written as
$$i\int_{\partial Q} a_j - i\int_{g^{-1}z}^z a_j,$$
where $Q$ denotes the geodesic quadrilateral with vertices at
$u$, $g^{-1}u$, $g^{-1}z$ and $z$.
By Stokes' Theorem the first integral can also be written as
$\int_Q da_j$, and this vanishes as we chose $a_j$ to be the harmonic
representative of its class, leaving just
$$i\int_{g^{-1}z}^z a_j.$$

Using this last equation
in the formula for the action of the imprimitivity algebra
(preceding Lemma 7) we see that we have a densely defined derivation
on the algebras $\ca$ and  $\ck$ because
the action can now be written as the
commutator of $\alpha$ with multiplication by the function
$\Omega_j(z) = i\int_u^z a_j$:
$$\delta_j\alpha = [\Omega_j,\alpha].$$
Suppose that $\alpha$ is a kernel
 {\it decaying rapidly}. By this we mean that it
 satisfies an estimate
$$
|\alpha (x,y)| \le \phi(d(x,y)),
$$
where $\phi$ is a positive and rapidly decreasing function on $\mathbb R$.
We claim that $\delta_j\alpha$ lies in $\ca$ or $\ck$ respectively: this
follows by noting that the map from $\hyp$ to
 $\real^{2g}$ given by $\Xi:z\mapsto (\Omega_1(z),\ldots,\Omega_{2g}(z))$
is the lift to $\hyp$ of the Jacobi map,
 \cite{GH}
(this map is usually regarded as mapping from $\Sigma_g$ to the Jacobi variety
$J(\Sigma_g)$, however we are thinking of it as a map between the universal
covers of these spaces). Now $\Xi$ is globally Lipschitz and this means
we may estimate the kernel of
$ [\Omega_j,\alpha]$. A simple argument shows that it
 also decays rapidly.

That this commutator also has the correct properties to define a derivation
on $\ca^\Gamma$ or $\ck^\Gamma$ follows from the
 fact that for $\gamma\in \Gamma$,
$$\Omega_j(\gamma.z)-\Omega_j(z)$$ is constant independent of $z$
so that $\Gamma$ equivariance or invariance is preserved.
In the case of the torus thought of as a rectangle in $\real^2$ with opposite
edges identified, one may take $a_j+ia_k = dz$, and then
$$(\delta_j+i\delta_k)\alpha = i[z-u,\alpha].$$
Thus our argument reproduces the standard Kubo
formula \cite{Xia} in the Euclidean case.

Notice that our map $\Xi$ from $\hyp$ to $\real^{2g}$
gives the period lattice in $\real^{2g}$ (that is the lattice determined
by the periods of the harmonic forms $a_j$) to be  the standard
integer lattice $\ints^{2g}$ so that $J(\Sigma_g)=\real^{2g}/\ints^{2g}$.
We may summarise the previous discussion
as

\begin{lemma} For operators $A_0, A_1, A_2$ in
 $\ck^\Gamma$ whose integral kernels
are rapidly decaying we have cyclic cocycles defined by
$$c_{j,k}(A_0, A_1, A_2)=
\tr_{\ck^\Gamma}(A_0[\delta_jA_1,\delta_kA_2])
= \tr_{\ck^\Gamma}(A_0[\Omega_j, A_1][\Omega_k,A_2])$$
for $j,k=1,\ldots,2g$.
\end{lemma}

Each of these formulae for $c_{jk}$ could in principle be regarded
as giving a Kubo formula so that we appear to have an embarrassment
of riches. However each on their own cannnot be related to
the Chern character of the Fredholm module of the next section.
A clue as to what is happening is provided by noting that each
two form $a_j\wedge a_{j+g}$ is harmonic and hence is a multiple
of the area two form on $\Sigma_g$, thus there is certainly
some degeneracy here and we resolve it at the end of the next section.

\section{A Fredholm module}

We shall now assume that $X$ has a spin structure, and we write $\cs$ for the
spin bundle.
The representation of $\cag$ can then be extended to an
action on $\chh_\ca^\chi\otimes\cs$.
This module can be equipped with Fredholm structure by taking $F$ to be
Clifford multiplication by a suitable unit vector (to be explained below), and
using the product of the trace on $\chh_\ca^\chi$ and the graded trace on
the Clifford algebra.
(If $\varepsilon$ denotes the grading operator on the spinors then the graded
trace is just $\tr\circ\varepsilon$.)

The same module can also be described more explicitly: it splits into
$\chh_\ca^\chi\otimes\cs^+\oplus \chh_\ca^\chi\otimes\cs^-$ (with the
 superscripted sign indicating the eigenvalue of $\varepsilon$), and
this may be written as
$\chh_\ca^{\theta_1}\oplus\chh_\ca^{\theta_2}$.
The involution $F$ is then a matrix multiplication operator of the form
$$F = \left(\begin{array}{cc}
0 &f_1\cr f_2 &0\cr
\end{array}\right),$$
with $(f_j.\psi)(z,g) = f_j(z,g)\psi(z,g)$
for some suitable functions $f_j \in C_c(G\times G)$,
satisfying $f_1=f_2^{-1}$.
For consistency, we require that for any $\psi \in \chh_\ca^{\theta_1}$,
$f_1.\psi \in \chh_\ca^{\theta_2}$.
Since
$$(f_1.\psi)(zk,g) = f_1(zk,g)\psi(zk,g) = \theta_1(k)f_1(zk,g)\psi(z,g),$$
we demand that $\theta_2(k)f_1(z,g) = \theta_1(k)f_1(zk,g)$, or
$$f_1(zk,g) = (\theta_1^{-1}\theta_2)(k)f_1(z,g),$$
and $\theta_1^{-1}\theta_2$ is known directly from the structure of $\cs$.
(When $X$ is the hyperbolic plane it is the complex character describing the
action of $K$ on the complex tangent space to $X$ at $u$.)
A short calculation shows that
$$([f_j,\alpha].\psi)(z,g) = \int_G
\left(f_j(z,g)-f_j(x^{-1}z,x^{-1}g)\right)
\frac{\alpha(zK,x)\psi(x^{-1}z,x^{-1}g)}
{\sigma(x,x^{-1}g)}\,dx.$$

We observe in the next subsection that this module is 2-summable
at least for kernels which decay sufficiently rapidly.  Assuming this fact
then it follows that
$((\omega*[f_j,\alpha]*[f_k,\beta]).\psi)(z,g)$ is given by
\begin{eqnarray*}
\int_{G\times G \times G}
\left(f_j(x^{-1}z,x^{-1}g)-f_j(y^{-1}z,y^{-1}g)\right)
\left(f_k(y^{-1}z,y^{-1}g)-f_k(u^{-1}z,u^{-1}g)\right)\\
\frac{\omega(zK,x)\alpha(x^{-1}zK,x^{-1}y)\beta(y^{-1}zK,y^{-1}u)}
{\sigma(x,x^{-1}g)\sigma(x^{-1}y,y^{-1}g)\sigma(y^{-1}u,u^{-1}g)}
\psi(u^{-1}z,u^{-1}g)\,dx\,dy\,du.
\end{eqnarray*} From
 this (and using $\tr$ to denote the usual trace on operators on
our module) we can calculate the cyclic cocycle on $\ca$ as
$$\tau_{c}(\omega,\alpha,\beta)
= \tr\left[\epsilon\omega*[F,\alpha]*[F,\beta]\right]
= \tr\left[\omega*
\left([f_1,\alpha]*[f_2,\beta]-[f_2,\alpha]*[f_1,\beta]\right)\right],$$
which can be expressed as
$$\int_{G/K\times G\times G}
\Phi(zK,x,y)\frac{\omega(zK,x)\alpha(x^{-1}zK,x^{-1}y)\beta(y^{-1}zK,y^{-1})}
{\sigma(x,x^{-1}y)\sigma(y,y^{-1})}\,dzK\,dx\,dy,$$
where
$$\Phi(zK,x,y) = \int_G
\left(f_1(x^{-1}z,x^{-1}g)-f_1(y^{-1}z,y^{-1}g)\right)
\left(f_2(y^{-1}z,y^{-1}g)-f_2(z,g)\right)$$
$$ \phantom{aaaaaaaaaa}- \left(f_2(x^{-1}z,x^{-1}g)-f_2(y^{-1}z,y^{-1}g)\right)
\left(f_1(y^{-1}z,y^{-1}g)-f_1(z,g)\right)\,dg.$$
(Using the equivariance of $f_j$ it is easy to check that this depends on
$z$ only through $zK$.)
Simplifying and using $f_1f_2 = 1$, the integrand reduces to
$$
\left(f_1(z,g)f_2(x^{-1}z,x^{-1}g)+f_1(x^{-1}z,x^{-1}g)f_2(y^{-1}z,y^{-1}g)
+f_1(y^{-1}z,y^{-1}g)f_2(z,g)\right)$$
$$-\left(f_2(z,g)f_1(x^{-1}z,x^{-1}g)+f_2(x^{-1}z,x^{-1}g)f_1(y^{-1}z,y^{-1}g)
+f_2(y^{-1}z,y^{-1}g)f_1(z,g)\right).$$
This can also be written more compactly as
$$\left|
\begin{array}{ccc}
1 &1 &1\\
 f_1(z,g) &f_1(x^{-1}z,x^{-1}g) &f_1(y^{-1}z,y^{-1}g)\\
 f_2(z,g) &f_2(x^{-1}z,x^{-1}g) &f_2(y^{-1}z,y^{-1}g)\\
\end{array}\right|,$$
or as
$$(1-f_1(x^{-1}z,x^{-1}g)f_2(z,g))(1-f_1(y^{-1}z,y^{-1}g)f_2(x^{-1}z,x^{-1}g))
(1-f_1(z,g)f_2(y^{-1}z,y^{-1}g)),$$
which also arises naturally from an alternative expression for
the cocycle.

Suppose that $\varphi$ is a $U(1)$ valued function on the group, which satisfies
$\varphi(kgh) = \chi_1(k)\varphi(g)\chi_2(h)$ for $k$ and $h$ in $K$ and some
$\sigma$-characters $\chi_1$ and $\chi_2$ of $K$.
If $\chi_1^{-1}\chi_2 = \theta_1^{-1}\theta_2$ we may take
$f_1(z,g) = \varphi(z^{-1}gz)$ to obtain a function satisfying
our earlier consistency condition.

In the case of $G = \su$ and $K$ the diagonal subgroup, we may take
the function $\varphi$ used by Connes \cite{Co2},
which is essentially the Mishchenko element.
With the group elements all conjugated
by $z$ it now follows as in \cite{Co}
 that $\Phi(zK,x,y)/4\pi i$ is the area of the hyperbolic geodesic
triangle with vertices $u$, $z^{-1}yz.u$ and $z^{-1}xz.u$.
Acting with $z$  and recalling that, since $u$ is stabilised by $K$,
$z.u$ can be identified with $s=zK$, $\Phi(zK,x,y)/4\pi i$ is
also the area of the
geodesic triangle with  vertices $s$, $y.s$ and $x.s$.

In the next subsection we will see that the module is 2-summable
for suitably decaying kernels.
Since
$f_j(\gamma^{-1}z,\gamma^{-1}g\gamma) =
\varphi(z^{-1}\gamma\gamma^{-1}g\gamma\gamma^{-1}z) = f_j(z,g)$,
$F$ preserves the $\Gamma$-invariant subspace, so that there
is a similar expression for a cyclic cocycle $\tau_{c,\Gamma}$ in
that case, except that $s$ is integrated only over the $\Gamma$ orbits in
$G/K$. More precisely, using $\tr_\Gamma$ to denote this
restricted range of integration, one has:

\begin{thm}
There is a 2-summable Fredholm module $(F, \chh_\tau^\chi\otimes\cs)$
over a dense subalgebra  ${\mathcal A}^\Gamma_0$ of ${\mathcal A}^\Gamma$,
stable under the holomorphic functional calculus, whose Chern character
is given by the area cocycle on $\mathbb H$. That is, in the
notation above, one has
$$
\tau_{c,\Gamma}(\omega, \alpha, \beta) =
-\tr_{\Gamma}\left[\epsilon\omega*[F,\alpha]*[F,\beta]\right]
= -\tr_{\Gamma}\left[\omega*
\left([f_1,\alpha]*[f_2,\beta]-[f_2,\alpha]*[f_1,\beta]\right)\right]
$$
which can be expressed as
$$
-\int_{\Gamma\backslash G/K\times G\times G}
\Phi(\Gamma
zK,x,y)\frac{\omega(zK,x)\alpha(x^{-1}zK,x^{-1}y)\beta(y^{-1}zK,y^{-1})}
{\sigma(x,x^{-1}y)\sigma(y,y^{-1})}\,dzK\,dx\,dy,
$$
where $\Phi$ is given as above. Therefore by the
index pairing in \cite{Co2}, one has
$$
 \Index(PFP) =    \langle[\tau_{c,\Gamma}], [P]\rangle ,
$$
where $P$ denotes a projection in  ${\mathcal A}^\Gamma_0$
and $\Index(PFP)$ denotes the index of the Fredholm operator $PFP$
acting on the Hilbert space $P\chh_\tau^\chi\otimes\cs$.
\end{thm}

We will prove Theorem 3 in the next subsection.
The version of Theorem 3 which applies to $\ck^\Gamma$
is as follows:

\begin{thm} There is a dense subalgebra ${\mathcal B}^\Gamma_0$
of  ${\mathcal B}^\Gamma$ stable under the holomorphic
functional calculus and a 2-summable Fredholm module
$(F, \chh_\tau^\chi\otimes\cs)$
for ${\mathcal B}^\Gamma_0$ with corresponding cyclic 2-cocycle
$$\tau_{c,\Gamma}(T_\omega,T_\alpha,T_\beta)
=-\int_{X_\Gamma\times X\times X}
\Phi(z,x,y)\varpi(z,x,y)T_\omega(z,x)T_\alpha(x,y)T_\beta(y,z)\,dz\,dx\,dy.$$
The  character of this Fredholm module, for
$P$ a projection in ${\mathcal B}^\Gamma_0$,
is given in the notation of Corollary 12 of Section 10, by
$$
   \Index(PFP) = 2(g-1)(\mbox{rank}\ \mathcal{E}^0 -
      \mbox{rank}\ \mathcal{E}^1)\in\mathbb{Z},
$$
where $\Index(PFP)$ again denotes the index of the Fredholm operator $PFP$
acting on the Hilbert space $P\chh_\tau^\chi\otimes\cs$.
\end{thm}

This theorem can be interpreted as
an {\em index theorem} equating an
analytic index with a topological index.
 Theorem
4 may be used to obtain the following result:

\begin{cor} Let $P$ be a projection into a gap in the spectrum
of the Hamiltonian $H_{\eta, V}$. Then
 $P$
 lies in a 2-summable dense subalgebra ${\mathcal B}^\Gamma_0$
of ${\mathcal B}^\Gamma$ so that
 in the notation of Corollary 12 of Section 10, one has
\begin{eqnarray*}
   \Index(PFP)  & = &  \langle\tau_{c,\Gamma}, [P]\rangle\\
& = & 2(g-1)(\mbox{rank}\ \mathcal{E}^0 -
      \mbox{rank}\ \mathcal{E}^1)\in\mathbb{Z}.
\end{eqnarray*}
\end{cor}

The statements referring to
Section 10 will be clear after we
establish there the hyperbolic analogues of
Xia's results \cite{Xia}. The proof of the  claim that the
spectral projections (corresponding to gaps in the
spectrum) of the Hamiltonian lie in  ${\mathcal B}^\Gamma_0$
and the  proof of Theorem 3 are
 contained in the next subsection.

\subsection{Proof of summability of the Fredholm module}

Here we discuss the technicalities needed for the proof of theorem 3
and of Corollary 6.
It is easy to calculate from the formulae in the previous
section that
2-summability requires finiteness of the expression
${\tr_\Gamma([f_1,\alpha]^*[f_2,\alpha])}$, which equals
$$\int_{X_\Gamma\times G\times G}[f_1(z,g) - f_1(x^{-1}z,x^{-1}g)]
[f_2(z,g) - f_2(x^{-1}z,x^{-1}g)]|\alpha(x^{-1}z,x^{-1})|^2\,dz\,dx\,dg,$$
where $z\in X_\Gamma$
 is a fundamental domain in $\hyp$ for the $\Gamma$ action.
Letting $p:G/K\rightarrow G$ be a cross section and
using the formulae for $f_1$ and $f_2$ this reduces to
$$\int_{X_\Gamma\times G\times G}
|\varphi(p(z)^{-1}gz)-\varphi(p(z)^{-1}gx^{-1}z)|^2|\alpha(x^{-1}z,x^{-1})|^
2\,dz\,dx\,dg.$$

We are more interested in the algebra of twisted kernels
as in Theorem 4, so we will present the argument
for them noting that the relation
for $\alpha$ in terms of $T_\alpha$ as given in the proof of Lemma 7
gives, by the the unitarity of
$\tau$ and $\varphi$,
$$|\alpha(x^{-1}z,x^{-1})|^2 = |T_\alpha(x^{-1}z,z)|^2.$$
Thus the summability result for $\ck^\Gamma$ implies that for $\ca^\Gamma$.

Making this substitution we then get for our integral
$$\int_{X_\Gamma\times G\times G}
|\varphi(p(z)^{-1}gz)-\varphi(p(z)^{-1}gx^{-1}z)|^2|T_\alpha(x^{-1}z,z)|^2\,
dz\,dx\,dg.$$
Finally we note that setting $p(z)^{-1}g = v^{-1}$ and
 $x^{-1}z = w$ it is clear that
$|\varphi(v^{-1}z)-\varphi(v^{-1}w)|$ depends only on the cosets $vK$ and $wK$,
allowing us to reduce the integral to
$$\int_{X_\Gamma\times X\times X}
|\varphi(v^{-1}z)-\varphi(v^{-1}w)|^2|T_\alpha(w,z)|^2\,dz\,dwK\,dvK.$$
Only the first factor depends on $v$.
Write  $v=\gamma v_0$ for $\gamma\in \Gamma$ and
$v_0\in X_\Gamma$. Then
we obtain for our integral after a change of variables:
$$\sum_{\gamma\in\Gamma}\int_{\gamma.X_\Gamma\times X\times X_\Gamma}
|\varphi(v_0^{-1}z)-\varphi(v_0^{-1}w)|^2|T_\alpha(\gamma w,\gamma
z)|^2\,dz\,dwK\,dv_0K,
$$
so that by the $\gamma$ invariance of the kernel $T_\alpha$
one obtains
$$\int_{X\times X\times X_\Gamma}
|\varphi(v_0^{-1}z)-\varphi(v_0^{-1}w)|^2|T_\alpha(w,z)|^2\,dz\,dwK\,dv_0K.$$
By a further change of variable we obtain
$$\int_{X\times X\times X_\Gamma}
|\varphi(z)-\varphi(w)|^2|T_\alpha(v_0^{-1}w,v_0^{-1}z)|^2\,dz\,dwK\,dv_0K.$$
Notice that, by lemma 6, $T_\alpha(v_0^{-1}w,v_0^{-1}z)$
is the integral kernel for the operator obtained by conjugating
by $U(v_0)$. It follows therefore that finiteness of the triple integral is
guaranteed by the convergence of
$$\int_{X\times X}
|\varphi(z)-\varphi(w)|^2|T_\alpha(w,z)|^2\,dz\,dwK.$$

To avoid repetition let us first focus on the case of greatest interest
where we consider the integral kernels of a spectral
projection $P$ of the Hamiltonian $H+V$
corresponding to a gap in the spectrum. As it is obtained
from the Hamiltonian using the smooth functional calculus
from a function of compact support we can
 obtain a growth estimate on the integral kernel
(see below) which will ensure convergence.
To lighten the notation we let $z,w\in X$ and
$k(z,w)$ denote the integral kernel as a function on
$X\times X$.
Taking $\chi_1 = 1$ so that $\psi(gk) = \varphi(g)$ is well defined,
the discussion of the previous paragraph
leads us to consider whether
$$\int\!\!\int |(\psi(z)-\psi(w))k(z,w)|^2 dz\,dw \eqno(*)$$
is finite.
Let $X_0$ be a fundamental domain for the diagonal action of $\Gamma$ on
$X \times  X$.    Then the previous integral, for $\Gamma$ invariant
kernels, is given by
$$\sum_{\gamma\in \Gamma}
\int\!\!\int_{X_0}|(\psi(\gamma.z)-\psi(\gamma.w))k(z,w)|^2dz\,dw.
$$
By an argument due to Connes \cite{Co}, we have the estimate
$$|\psi(\gamma.z)-\psi(\gamma.w)|^2 \leq C\exp(-2d(u,\gamma.z)+C_1 d(\gamma.z,
\gamma.w)).$$
(Here $u$ is the base point in $X$, $d$  denotes the hyperbolic metric
and $C$,$C_1$ are constants.)

We claim that, in addition, the following estimate holds:
$$ |k(z,w)|^2\leq C_2\exp(-C_3d(z,w)^2), \eqno(**)$$
where $C_2,C_3$ are constants.
This fact goes back to \cite{CGT}
although in the form we need it here, for the Hamiltonians $H+V$
of Section 4, it can be deduced from \cite{BrSu}. This is because
\cite{BrSu} prove (**) when $k$ is the kernel of the heat operator
$e^{-(H+V)}$. Now by the argument of Corollary 3 there is a smooth
function of compact support $g$ such that $g(e^{-(H+V)})$ is the
spectral projection $P$. To prove that the kernel of $g(e^{-(H+V)})$
satisfies (**) it suffices to observe that we can approximate
$g$ uniformly by polynomials without constant term so that the kernel of
 $g(e^{-(H+V)})$ has the same off-diagonal decay estimate as
the kernel of $e^{-(H+V)}$, namely (**).

Hence the integral in (*) above is smaller than
$$\sum_{\gamma\in\Gamma}\int\!\!\int_{X_0} C_4 \exp(-C_3d(z,w)^2+C_1d(z,w)-
2d(u,\gamma.z))\eqno(***)$$
for suitable constants $C_j, j=1,2,3$.
As the area in hyperbolic space grows like $\exp(d(u,z))$,
convergence of the infinite sum in (***) is handled by the convergence of the
Poincare series
$\sum_{\gamma\in\Gamma} \exp(-2d(u,\gamma.z))$.
The convergence of the integral in (***), over the fundamental domain,
is handled by the
exponential factor involving the square of the hyperbolic distance
and noting that the integration in the diagonal direction in $X_0$ is
over a finite range.
(It is also possible to prove 2-summability in the case when $\Gamma$ is
trivial by exploiting the fact that in that case one may use
kernels with restrictions on their support.)

Since operators with kernels which have
 support in a band around the diagonal are
dense in the algebras $\ca^\Gamma$ and $\ck^\Gamma$
so too is the set of operators with kernels satisfying (**).
Now the  finiteness of (*) is equivalent to asserting that $[F,T_\alpha]$
is Hilbert-Schmidt.

\noindent{\bf Definition}.
We denote by ${\mathcal B}^\Gamma_0$
the subalgebra
consisting of operators $A\in {\mathcal B}^\Gamma$,
with $[F,A]$ a Hilbert-Schmidt operator.

The argument of the previous paragraph shows that ${\mathcal B}^\Gamma_0$
is dense.
Now
by \cite{Co} ${\mathcal B}^\Gamma_0$ is stable under
 the holomorphic functional calculus.
A similar remark handles the existence of the analogous
dense subalgebra  ${\mathcal A}^\Gamma_0$ of ${\mathcal A}^\Gamma$.
This completes the proof of Theorem 3 and the claim
concerning the spectral projections of the Hamiltonian as we promised.

\subsection{The hyperbolic Connes-Kubo formula}

Our aim in this subsection is to give a geometric interpretation
to the cocycles defined in lemma 12
and to prove that a suitable linear combination of them
is cohomologous to the cocycle $\tau_{c,\Gamma}$ arising from the
Fredholm module $(F, \chh_\tau^\chi\otimes\cs)$.

 To do this we begin
by introducing, for operators $A_0, A_1, A_2$ in $\ck^\Gamma$
whose kernels $k_0,k_1,k_2$ are exponentially decaying
 (cf equation (**) of the previous
subsection) the cyclic cocycle $c_K$ defined by
\begin{eqnarray*}
c_K(A_0, A_1, A_2) & = & \sum_{j=1}^g c_{j,j+g}(A_0, A_1, A_2)\\
 & = & \sum_{j=1}^{g}\int_{X_\Gamma\times X\times X}\varpi(z,x,y)
\Psi_j(z,x,y)k_0(z,x)k_1(x,y)k_2(y,z)\,dz\,dx\,dy,
\end{eqnarray*}
where
$$\Psi_j(z,x,y)=(\Omega_j(x)-\Omega_{j}(y))(\Omega_{j+g}(y)-\Omega_{j+g}(z))-
(\Omega_{j+g}(x)-\Omega_{j+g}(y))(\Omega_{j}(y)-\Omega_{j}(z)).$$
We claim that
$\sum_{j=1}^g\Psi_j(z,x,y)$ is proportional to the `symplectic area'
of a triangle in $\real^{2g}$ with vertices $\Xi(x),\Xi(y),\Xi(z)$.

To prove this it suffices to assume that the origin is one of the vertices of
the triangle, so suppose $z$ is the base point in
$\hyp$. Then we need to consider the expression
$$\sum_{j=1}^g\Psi_j(z,x,y)=\sum_{j=1}^g(\Omega_j(x)\Omega_{j+g}(y)-\Omega_{
j+g}(x)\Omega_{j}(y)).$$
Let $s$ denote the symplectic form on $\real^{2g}$ given by:
$$s(u,v)=\sum_{j=1}^g(u_jv_{j+g} -u_{j+g}v_j).$$
The so-called `symplectic area' of a triangle
with vertices $0,\Xi(x),\Xi(y)$ may be seen to be
$s(\Xi(x),\Xi(y))$. To appreciate this, however,
we need to utilise an argument from
{\cite{GH} (pp 333-336)}.
In terms of the standard basis of $\real^{2g}$ (given in this case
by vertices in the integer period lattice arising from our choice of basis
of harmonic one forms) and corresponding coordinates $u_1,u_2,\ldots u_{2g}$
the form $s$ is the two form on $\real^{2g}$ given by
$$\omega_J=\sum_{j=1}^g du_j\wedge du_{j+g}.$$
Now the  `symplectic area' of a triangle in $\real^{2g}$ with
 vertices $0,\Xi(x),\Xi(y)$
is given by integrating $\omega_J$ over the triangle and
a brief calculation reveals that this yields
$s(\Xi(x),\Xi(y))/2$, proving our claim.

The previous argument establishes the following result.

\begin{prop} The higher genus analogue of the Kubo formula is given by
the cyclic
cocycle $\tau_K$ on $\ck^\Gamma$ defined by
\begin{eqnarray*}
\tau_K(A_0, A_1, A_2) & = & \sum_{j=1}^g c_{j,j+g}(A_0, A_1, A_2)\\
 & = & \sum_{j=1}^{g}\int_{X_\Gamma\times X\times X}
\Psi_j(z,x,y)\varpi(z,x,y)k_0(z,x)k_1(x,y)k_2(y,z)\,dz\,dx\,dy.
\end{eqnarray*}
Here the $k_j$ are the kernels of the $A_j, j=0,1,2$
(three exponentially decaying elements of $\ck^\Gamma$) and
$\sum_{j=1}^{g}\Psi_j(z,x,y)$ is proportional to the
`symplectic area' of the Euclidean triangle
$\Delta_E$ in $\real^{2g}$ with vertices
 $\Xi(x),\Xi(y),\Xi(z)$.
\end{prop}

To compare the cocycle $\tau_K$ with the cocycle $\tau_{c,\Gamma}$
arising from our Fredholm module
we note that
 the pull back form $\Xi^*(du_j)$ is dual to the homology cycle $A_j$ for
$j=1,\ldots, g$ and dual to $B_{j-g}$ for $ j=g+1,\ldots, 2g$ (cf\cite{GH}).
Thus $\Xi^*(du_j)$ differs from $a_j$ by an exact one form.
Hence
 $\Xi^*(\omega_J)$ differs from $\sum_{j=1}^g a_j\wedge a_{j+g}$
by an exact two form. But each term $a_j\wedge a_{j+g}$
is harmonic and hence proportional to the two form $\omega_\hyp$
on $\hyp$. So we have for some constant $\kappa$,
and geodesic triangle $\Delta\subset\hyp$,
$$\int_\Delta\omega_\hyp=\kappa\int_\Delta\Xi^*(\omega_J) =
\kappa\int_{\Xi(\Delta)}\omega_J.
$$
Actually a calculation reveals that one can do a little better than this
and proves that
$$\kappa\Xi^*(\omega_J)=\omega_\hyp.$$
Now $\Xi$ cannot map geodesic triangles
to Euclidean triangles in $\real^{2g}$
as $\Xi({\Delta})$ is a compact subset of a non-flat embedded two
dimensional surface in $\real^{2g}$.
Moreover as  $\Psi_j(z,x,y)=0$ whenever the images of $z,x,y$ under $\Xi$
lie in a Lagrangian subspace (with respect to the symplectic
form $s$) of $\real^{2g}$, $\tau_K$ and $\tau_{c,\Gamma}$
are not obviously proportional.

After suitable normalisation we will, however, prove they are cohomologous.
First renormalise $\omega_J$ so that
$\Xi^*(\omega_J)=\omega_\hyp$
and then normalise $\sum_{j=1}^g\Psi_j(z,x,y)$ so that it equals
$-4\pi i\int_{\Delta_E}\omega_J$.
Next we write $\omega_J= d\theta$.
Considering the difference $\tau_K-\tau_{c,\Gamma}$ one sees that the
key is to understand
$$\int_{\Xi(\Delta)}\omega_J\quad -\int_{\Delta_E}\omega_J
=\int_{\partial\Xi(\Delta)}
\theta \quad - \int_{\partial\Delta_E}\theta.
$$
Now this difference of integrals around the boundary can be written as
the sum of three terms corresponding to splitting the boundaries
$\partial\Xi(\Delta)$ and $\partial\Delta_E$ into
three arc segments each. We introduce some notation for this,
writing
$$\partial\Xi(\Delta)=\Xi(\ell(x,y))\cup\Xi(\ell(y,z))\cup\Xi(\ell(z,x)),$$
where $\ell(x,y)$ is the geodesic in $\hyp$ joining $x$ and $y$
(with the obvious similar definition of the other terms). We also write
$$\partial\Delta_E=m(x,y)\cup m(y,z)\cup m(z,x),$$
where $m(x,y)$ is the straight line joining $\Xi(x)$ and $\Xi(y)$
(and again the obvious definition of the other terms).
Then we have
$$\int_{\partial\Xi(\Delta)}\theta\quad - \int_{\partial\Delta_E}\theta=
h(x,y)+h(y,z)+h(z,x)$$
where $h(x,y)= \int_{\Xi(\ell(x,y)}\theta - \int_{m(x,y)}\theta$
with similar definitions  for $h(y,z)$ and $h(z,x)$.

Notice that we can write $h(x,y)=\int_{D_{xy}}\omega_J$
where $D_{xy}$ is a disc with boundary $m(x,y)\cup \ell(x,y)$. From
 this it is easy to see that $h(\gamma x,\gamma y)=h(x,y)$
for $\gamma\in \Gamma$.
Introduce the bilinear functional $\tau_1$ on $\ck^\Gamma$
given by
$$\tau_1(A_0,A_1) =-4\pi i \int_{X_\Gamma \times X}
h(x,y)k_0(x,y)k_1(y,x)\,dx\,dy
=-4\pi i\tr_{\ck^\Gamma}(A_hA_1),$$
where, if $A_0$ has kernel $k_0(x,y)$, $A_h$ is the operator
with kernel $h(x,y)k_0(x,y)$.
Of course this definition begs the question of whether the trace is finite.

In order to prove that $\tau_1$ is densely defined
we start with some preliminary observations. By \cite{M+R+W} and \cite{Ren2}
there is an isomorphism
$$
\Phi_F: {\mathcal B}^\Gamma\cong C^*_r(\Gamma_g, \sigma)\otimes {\mathcal K}
(L^2(F)).
$$
Here $F$ denotes a fundamental domain for the action of $\Gamma_g$
on $\mathbb H$. (Note that by the Packer-Raeburn stabilization theorem, one has
$C^*_r(\Gamma_g, \sigma)\otimes {\mathcal K} \cong {\mathcal K} \rtimes_r
\Gamma_g$.) Now any element $x$ in $C^*_r(\Gamma_g, \sigma)\otimes {\mathcal K}$
can be written as a matrix $(x_{ij})$, where $x_{ij} \in C^*_r(\Gamma_g,
\sigma)$.
So we can define
$$N_k(x) = ( \sum_{i,j}\nu(x_{ij})^2)^{\frac{1}{2}},$$
where $$\nu(x_{ij}) = (\sum_{h\in \Gamma_g} (1+\ell(h)^{2k})
|x(h)|^2)^{\frac{1}{2}}$$
and $\ell$ denotes the word length function on the group $\Gamma_g$.
By a mild modification of the argument given in \cite{Co2},
III.5.$\gamma$, one can prove that there is a subalgebra
${\mathcal B}^\Gamma_\infty$ of ${\mathcal B}^\Gamma$ which contains
${\mathbb C}(\Gamma_g, \sigma)\otimes {\mathcal R}$, where $\mathcal R$
denotes the algebra of smoothing operators on $F$, is stable
under the holomorphic functional calculus, and is such that $N_k(x)<\infty$
for all $x\in {\mathcal B}^\Gamma_\infty$ and $k\in \mathbb N$. Then one
shows as
in \cite{Co2} that the trace $\tau \otimes \mbox{Tr}$ on
${\mathbb C}(\Gamma_g, \sigma)\otimes {\mathcal R}$, is continuous
for the norm $N_k$, for $k$ sufficiently large, and thus extends by continuity
to ${\mathcal B}^\Gamma_\infty$. Note that elements in
${\mathcal B}^\Gamma_\infty$ have Schwartz kernels which have
rapid decay away from the diagonal.
An alternate equivalent construction of
${\mathcal B}^\Gamma_\infty$ would be to use the algebra $A_{g,\sigma}$ as
in Section 10, and the results of \cite{Ji}.

Summarizing this, we have

\begin{prop} The algebra ${\mathcal B}^\Gamma_\infty$ is dense
in ${\mathcal B}^\Gamma$, is closed under the
holomorphic functional calculus
and is contained in the ideal ${\mathcal I}$ of
$\ck^\Gamma$ consisting of operators with finite trace.
\end{prop}

Now $\tau_K$ is defined on  ${\mathcal B}^\Gamma_\infty$ while
$\tau_{c,\Gamma}$ is defined on  ${\mathcal B}^\Gamma_0$
as we noted earlier. Both of these algebras contain the
operators whose Schwartz kernels are supported in a band around the
diagonal. Thus the subalgebra
${\mathcal B}^\Gamma_\infty\cap {\mathcal B}^\Gamma_0$ is dense
and stable under the holomorphic functional calculus.
If $b$ denotes the Hochschild coboundary map then a straightforward
calculation reveals that $b\tau_1= \tau_K-\tau_{c,\Gamma}$.
The Lipschitz property of the Jacobi map means that $h(x,y)$
grows at worst like the square of the hyperbolic distance from $x$ to $y$
so that if
$A_0\in\ck^\Gamma_\infty$ then so too does $A_h$.
Hence we have $\tau_1$ defined on $\ck^\Gamma_\infty\cap{\mathcal B}^\Gamma_0$
proving the following theorem.

\begin{thm}
The Kubo cocycle $\tau_K$ and the Chern character
cocycle $\tau_{c,\Gamma}$ arising
as the Chern class of the
Fredholm module $(F, \chh_\tau^\chi\otimes\cs)$, are cohomologous
as cyclic cocycles on $\ck^\Gamma_\infty\cap{\mathcal B}^\Gamma_0$.
\end{thm}

This theorem replaces the Connes-Kubo formula in genus one. The latter formula
states that the two cocycles of the theorem are equal. We see that the
 situation is more complex for genus $g$ but from the viewpoint of $K$-theory
as described in the next section this theorem is enough to
give integrality of the Hall conductivity defined either from the character of
our Fredholm module or from the hyperbolic Kubo formula.

\section{$K$-theory aspects}

In this section, we compute the $K$-groups of the twisted group
$C^*$-algebras which are relevant to the quantum Hall effect on the
hyperbolic plane
as a special case of more general theorems about the $K$-groups of the twisted
group $C^*$-algebras of groups $\Gamma$ which are uniform lattices
 in $K$-amenable
Lie groups.

We recall that any solvable Lie group, and in fact any amenable Lie group
is $K$-amenable.
However, it has been
proved by Kasparov \cite{Kas1} in the case of the
non-amenable groups ${\mathbf{SO}}_0(n,1)$ and by Julg-Kasparov
\cite{JuKas} in the case of ${\mathbf{SU}}(n,1)$
 that these are $K$-amenable Lie groups.
Cuntz \cite{Cu} has shown that the class of $K$-amenable groups
is closed under the operations of taking subgroups, under free products
and under direct products.
Our method uses the $K$-amenability results of Kasparov \cite{Kas1} and the
Packer-Raeburn stabilization theorem \cite{PR}.
In \cite{PR1}, one can find an example where the twisted $K$-theory
$K^{*}(\Gamma\backslash G/K, \delta(B_\sigma))$ is {\em not}
isomorphic to $K^{*}(\Gamma\backslash G/K)$, even when $G$ is the $K$-amenable
solvable group $\mathbb R^n \rtimes \mathbb R$, $K=\{e\}$,  $\Gamma=
\mathbb Z^n
\rtimes
\mathbb Z$ and for some multiplier $\sigma$ on $\Gamma$ with non-trivial
Dixmier-Douady
invariant $\delta(B_\sigma) = \delta(\sigma) \ne 0$. However,
$K^{*}(\Gamma\backslash G/K, \delta(B_\sigma))$ is
isomorphic to $K^{*}(\Gamma\backslash G/K)$
whenever the Dixmier-Douady invariant
$\delta(B_\sigma) = \delta(\sigma) = 0$ is
trivial. We identify the Dixmier-Douady
invariant $\delta(B_\sigma)$ with the image of $\sigma$ under the connecting
homomorphism $\delta : H^2(\Gamma, U(1)) \to H^3(\Gamma, \mathbb Z)$ of the
change of coefficients exact sequence in cohomology, corresponding to the
short exact sequence
of coefficient groups
\[
   1\to\mathbb{Z} \overset{i}{\to} \mathbb{R}
      \overset{e^{2\pi\sqrt{-1}}}{\longrightarrow} U(1) \to 1.
\]
This enables us to prove vanishing theorems for the Dixmier-Douady invariant
whenever $\Gamma$ is a lattice in a connected
Lie group $G$ such that $\dim(G/K) \le 3$, where $K$ is a maximal compact
subgroup
of $G$, and therefore we obtain in this case
$$
K_*(C^*(\Gamma,\sigma)) \cong K^{* + \dim(G/K)}(\Gamma\backslash G/K),
$$
where $\sigma$ is any multiplier on $\Gamma$.
This is the case for the Riemann surfaces which are the
object of our study in this paper.

We begin by reviewing the concept of $K$-amenable groups.
Let $G$ be a connected Lie group and $K$ be a maximal compact subgroup.
For our purposes, we will assume that $\dim(G/K)$ is even, and that it has
a $G$-invariant $\Spin^{\mathbb{C}}$ structure.  Using the
$\Spin^{\mathbb{C}}$ structure, we can form the $G$-invariant
Dirac operator $\slash\!\!\!\partial$ on $G/K$.  It is a first order,
self-adjoint,
elliptic differential operator acting on $L^2$ sections of the
$\mathbb{Z}_2$ graded homogeneous bundle of spinors $S$.  Consider $F =
\slash\!\!\!\partial \ (1+  \slash\!\!\!\partial^2)^{-1/2}$, which
is a $0$th order pseudo-differential operator acting on $H = L^2(G/K, S)$.
$C_0(G/K)$ acts on $H$
by multiplication operators, $f\to M_f$.  Also $G$ acts on $C_0(G/K)$ and
on $H$ by left translation, and $F$ is $G$-invariant.  Therefore $(H,M,F)$
defines a canonical element, called the \emph{Dirac element},
\[
   \alpha_G \in KK_G (C_0(G/K), \mathbb{C}).
\]

\begin{thm}[\cite{Kas2}]
There is a canonical element, called the \emph{Mishchenko element}
\[
   \beta_G \in KK_G(\mathbb{C}, C_0(G/K)),
\]
such that one has the following intersection products:
\begin{enumerate}
\item[(1)] $\alpha_G \otimes_{\mathbb{C}} \beta_G = 1_{C_0(G/K)} \in
KK_G(C_0(G/K), C_0(G/K))$
\item[(2)] $\beta_G \otimes_{C_0(G/K)} \alpha_G = \gamma_G \in
KK_G(\mathbb{C}, \mathbb{C})$ where $\gamma_G$ is an idempotent in
$KK_G(\mathbb{C}, \mathbb{C})$.
\end{enumerate}
\end{thm}

The Mishchenko element $\beta_G$ can be described as follows.  First assume
that either $G$ is semisimple or that $G=\mathbb{R}^n$.  Then the Killing
form on $G$ defines a $G$-invariant Riemannian metric of non-positive
sectional curvature on $G/K$.  Let $\mathcal{E} = C_0(G/K, \mathcal{S}^*)$
be the
space of continuous sections of the dual spin bundle $\mathcal{S}^*$ which
vanish
at infinity.  Let $F$ be a bounded operator on $\mathcal{E}$ defined as
\[
   F\xi (x) = c(V(x,x_0))\xi(x),
\]
where $\xi\in\mathcal{E},\ V(x,x_0)\in T_x(G/K)$ is the unit vector which
is tangent to the unique geodesic from $x_0 \in G/K$ to $x$ and $c(V(x,
x_0))$ denotes Clifford multiplication by $V(x,x_0)$.  Then $V(x,x_0)$ is
well defined outside a small neighbourhood of $x_0$ and can be extended
continuously in any way to all of $G/K$. As $F$ is adjointable it lies in
$\mathcal{L}(\mathcal{E})$.  Also since
\[
   F^2 \xi(x) = \| V(x,x_0) \|^2 \xi(x),
\]
we see that $F^2-1\in \mathcal{K}(\mathcal{E})$ is a compact operator in
$\mathcal{L}(\mathcal{E})$.  For $g\in G$, define
\[
   (g.F)\xi (x) = c(V(x, gx_0))\xi(x).
\]
Since $G/K$ has negative sectional curvature, the function on $G/K$ defined
by
\[
   x \to \| V(x,x_0) - V(x,x_1) \|,\quad x_0, x_1 \in G/K,
\]
vanishes at infinity, and so is in $C_0(G/K)$.
 Therefore $g.F - F \in \mathcal{K}(\mathcal{E})$ and $(\mathcal{E}, F)$
defines an
element $\beta_G \in KK_G (\mathbb{C}, C_0(G/K))$.
The Mishchenko element
$\beta_G$ is constructed by induction in the general case.

\begin{thm}[\cite{Kas2}]
If $G$ is amenable, then $\gamma_G = 1$.
\end{thm}

This motivates the following definition  (\cite{Kas2}).

\begin{defn*}
A Lie group $G$ is said to be $K$-amenable if $\gamma_G$ = 1.
\end{defn*}

\begin{thm}[\cite{Kas1}\cite{JuKas}]
The non-amenable groups ${\mathbf{SO}}(n, 1)$ and ${\mathbf{SU}}(n,1)$ are
$K$-amenable.
\end{thm}

Let $\Gamma\subset G$ be a lattice in $G$ and $A$ be an algebra
admitting an automorphic action of $\Gamma$.
Then the cross product algebra $[A\otimes C_0(G/K)]\rtimes \Gamma$,
is Morita equivalent to the algebra of continuous sections
vanishing at infinity
$C_0(\Gamma\backslash G/K, \mathcal{E})$, where $\mathcal{E}\to
\Gamma\backslash G/K$ is the flat $A$-bundle defined as the quotient
\begin{equation}
\tag{*}
   \mathcal{E} = (A\times G/K)/\Gamma \to \Gamma\backslash G/K.
\end{equation}
Here we consider the diagonal action of $\Gamma$ on $A\times G/K$.

\begin{thm}[\cite{Kas2}]
If $G$ is $K$-amenable, then
$(A\rtimes\Gamma)\otimes C_0(G/K)$ and $[A\otimes C_0(G/K)]\rtimes
\Gamma$
have the same $K$-theory.
\end{thm}

Combining Theorem 9 with the remarks above, one gets the following
important corollary.

\begin{cor}
If $G$ is $K$-amenable, then $(A\rtimes\Gamma)\otimes C_0(G/K)$ and
$C_0(\Gamma\backslash G/K, \mathcal{E})$ have the same $K$-theory.
Equivalently, one has for $j=0,1$,
$$
K_j(C_0(\Gamma\backslash G/K, \mathcal{E})) \cong K_{j+ \dim(G/K)}
(A\rtimes\Gamma).
$$
\end{cor}

We now come to the main theorem of this section,
which generalizes theorems of \cite{PR}, \cite{PR2}.

\begin{thm}
Suppose that $\Gamma$ is a lattice in a $K$-amenable Lie group $G$  and that
$K$ is a maximal compact subgroup of $G$. Then
$$
K_*(C^*(\Gamma,\sigma)) \cong K^{* + \dim(G/K)}(\Gamma\backslash G/K,
\delta(B_\sigma)),
$$
where $\sigma \in H^2(\Gamma, U(1))$ is any multiplier on $\Gamma$,
$K^{*}(\Gamma\backslash G/K,
\delta(B_\sigma))$ is the twisted $K$-theory of a continuous trace
 $C^*$-algebra
$B_\sigma$ with spectrum $\Gamma\backslash G/K$, and $\delta(B_\sigma)$
denotes the
Dixmier-Douady invariant of $B_\sigma$.\end{thm}

\begin{proof}
1. Taking the case $A=\mathbb{C}$ and the trivial action of $\Gamma$ on
$\mathbb{C}$, one sees by Corollary 7
that $C^*(\Gamma)$ and $C_0(\Gamma\backslash G/K)$
have the same $K$-theory when $\gamma_G=1$.

2. Let $\sigma\in H^2(\Gamma, U(1))$, then the twisted cross product
algebra $A\rtimes_\sigma \Gamma$ is stably equivalent to $(A\otimes
\mathcal{K})\rtimes\Gamma$ where $\mathcal{K}$ denotes compact operators.
This is the Packer-Raeburn stabilization trick \cite{PR}
(note that the $\Gamma$ action on $\mathcal{K}$ takes some time to describe
and we refer the reader to \cite{PR} for details).  Using Corollary 7 again,
one sees
that $A\rtimes_\sigma\Gamma \otimes C_0(G/K)$ and $C_0(\Gamma\backslash
G/K, \mathcal{E}_\sigma)$
have the same $K$-theory, whenever $G$ is $K$-amenable, where
$$
   \mathcal{E}_\sigma = (A\otimes\mathcal{K} \times G/K)/\Gamma \to
      \Gamma\backslash G/K
$$
is a flat $A\otimes\mathcal{K}$-bundle over $\Gamma\backslash G/K$ and
$K$ is a maximal compact subgroup of $G$.
In the particular case when $A= \mathbb{C}$, one sees that $C^*_r(\Gamma,
\sigma)\otimes C_0(G/K)$ and $C_0(\Gamma\backslash G/K,
\mathcal{E}_\sigma)$ have the same
$K$-theory whenever $G$ is $K$-amenable, where
\[
   \mathcal{E}_\sigma = (\mathcal{K} \times G/K)/\Gamma\to\Gamma\backslash G/K.
\]
But the twisted $K$-theory $K^{*}(\Gamma\backslash G/K, \delta(B_\sigma))$ is
by definition the $K$-theory of the continuous trace $C^*$-algebra $B_\sigma
= C_0(\Gamma\backslash G/K, \mathcal{E}_\sigma)$ with spectrum
$\Gamma\backslash G/K$.
Then
$$
K_*(C^*(\Gamma,\sigma)) \cong K^{* + \dim(G/K)}(\Gamma\backslash G/K,
\delta(B_\sigma)).
$$

\end{proof}

\begin{rems*}
Consider
 the flat case, when $G = \mathbb R^{2n} \rtimes \so(2n)$ is the Euclidean
group, $K= \so(2n)$,
and  $\Gamma \subset G$ is a Bieberbach group, that is, $\Gamma$ is a uniform
lattice in $G$. One can define a generalization of \lq\lq noncommutative
flat manifolds" by regarding
$C^*(\Gamma, \sigma)$ as such an object, where $\sigma$ is any group
2-cocycle on $\Gamma$,
by virtue of the fact that
$$
K_*(C^*(\Gamma,\sigma)) \cong K^{*}(\Gamma\backslash G/K).
$$
\end{rems*}

Our next main result says that for lattices in $K$-amenable Lie groups,
the reduced and unreduced twisted group $C^*$-algebras have canonically
isomorphic $K$-theories. Therefore all the results that we prove regarding
the $K$-theory of these reduced twisted group $C^*$-algebras are also valid
for the unreduced twisted group $C^*$-algebras.

\begin{thm}
Let $\sigma\in H^2(\Gamma, U(1))$ be a multiplier on $\Gamma$ and
$\Gamma$ be a lattice in a $K$-amenable Lie group. Then the canonical
morphism $C^*(\Gamma, \sigma) \rightarrow C^*_r(\Gamma, \sigma)$ induces
an isomorphism
$$
K_*(C^*(\Gamma, \sigma)) \cong K_*(C^*_r(\Gamma, \sigma)).
$$
\end{thm}

\begin{proof}
We note that by the Packer-Raeburn trick, one has
$$
C^*(\Gamma, \sigma) \otimes {\mathcal K} \cong  {\mathcal K}  \rtimes  \Gamma
$$
and
$$
C^*_r(\Gamma, \sigma) \otimes {\mathcal K} \cong  {\mathcal K}  \rtimes_r
\Gamma,
$$
where $\rtimes_r$ denotes the reduced crossed product.
Since $\Gamma$ is a lattice in a $K$-amenable Lie group,
the canonical
morphism ${\mathcal K}  \rtimes  \Gamma \rightarrow
{\mathcal K}  \rtimes_r  \Gamma$ induces
an isomorphism
$$
K_*({\mathcal K}  \rtimes  \Gamma) \cong K_*({\mathcal K}  \rtimes_r  \Gamma),
$$
which proves the result.
\end{proof}

We now specialize to the case when $G= \so_0(2,1)$, $K=\so(2)$ and
$\Gamma = \Gamma_g$ is the fundamental group of a Riemann surface of genus
$g>1$, \ $\Sigma_g$, where  $\Gamma_g \subset G$ and $G$ is $K$-amenable, or
when $G= \mathbb R^2$, $K= \{e\}$ and $g=1$, with $\Gamma_1$ being $\mathbb
Z^2$.

\begin{cor} 
Let $\sigma \in H^2(\Gamma_g, U(1))$ be any multiplier
on $\Gamma_g$.  Then
\begin{enumerate}
\item $K_0(C^*_r(\Gamma_g, \sigma)) \cong K_0(C^*_r(\Gamma_g)) \cong
K^0(\Sigma_g) \cong \mathbb{Z}^2$
\item $K_1(C_r^*(\Gamma_g, \sigma)) \cong K_1 (C_r^*(\Gamma_g)) \cong
K^1(\Sigma_g) \cong \mathbb{Z}^{2g}$.
\end{enumerate}
\end{cor}

\begin{proof}
In dimension 2 the Chern character is an isomorphism over the integers and
therefore we see that
\begin{align*}
   K^0 (\Sigma_g) &\cong H^0(\Sigma_g, \mathbb{Z}) \oplus H^2(\Sigma_g,
      \mathbb{Z}) \cong \mathbb{Z}^2, \\
\intertext{and that}
   K^1(\Sigma_g) &\cong H^1(\Sigma_g, \mathbb{Z}) \cong \mathbb{Z}^{2g}.
\end{align*}
By Theorem 10 we have
\[
   K_j(C^*_r(\Gamma_g)) \cong K^j(\Sigma_g)\quad\text{for }j=0,1,
\]
and
\[
   K_j(C_r^*(\Gamma_g,\sigma)) \cong K_j(\Sigma_g, \delta(B_\sigma)),\quad
      j = 0,1,
\]
where $B_\sigma = C(\Sigma_g, \mathcal{E}_\sigma)$.
Finally, because $\mathcal{E}_\sigma$ is a locally trivial flat bundle of
 $C^*$-algebras over $\Sigma_g$,
with fibre $\mathcal{K}$ ($=$ compact operators),
it has a Dixmier-Douady invariant $\delta(B_\sigma)$ which can be
viewed as the obstruction to $B_\sigma$ being Morita
equivalent to $C(\Sigma_g)$.  But
\[
   \delta(B_\sigma) = \delta(\sigma) \in H^3(\Sigma_g, \mathbb{Z}) = 0.
\]
Therefore $B_\sigma$ is Morita equivalent to $C(\Sigma_g)$
and we conclude that
\[
   K_j(C^*_r(\Gamma_g, \sigma)) \cong K^j(\Sigma_g)\quad j = 0,1.
\]
\end{proof}

\begin{cor} 
Let $G$ be a connected Lie group and $K$ a maximal compact subgroup such
that $\dim(G/K) = 3$.  Let $\Gamma$ be a uniform lattice in $G$ and
$\sigma\in H^2(\Gamma, U(1))$ be any multiplier on $\Gamma$.
If  $G$ is $K$-amenable, then
\begin{equation}
\tag{*}
   K_j(C^*_r(\Gamma,\sigma)) \cong K_j(C^*_r(\Gamma)) \cong
      K^{j+1}(\Gamma\backslash G/K),\qquad\mbox{for }j=0,1 \pmod 2.
\end{equation}
\end{cor}

\begin{proof}
By Theorem 10, we see that
\[
   K_j(C^*_r(\Gamma)) \cong K^{j+\dim(G/K)}(\Gamma\backslash G/K),
\qquad\mbox{for }j=0,1 \pmod 2.
\]
By the Packer-Raeburn stabilization trick, $C_r^*(\Gamma,\sigma)$ is Morita
equivalent to $\mathcal K \rtimes \Gamma$, and because $G$ is $K$-amenable,
$\mathcal K \rtimes \Gamma \otimes C_0(G/K)$ is Morita equivalent to
$B_\sigma = C(\Gamma\backslash G/K, \mathcal{E}_\sigma)$, where
$\mathcal{E}_\sigma$ is
as before, a locally trivial bundle of $C^*$-algebras over
$\Gamma\backslash G/K$
with fibre $\mathcal{K}$.  Finally, the
Dixmier-Douady invariant
\[
   \delta(B_\sigma) = \delta(\sigma) \in H^3(\Gamma\backslash G/K, \mathbb{Z})
      \cong H^3(\Gamma,\mathbb{Z}).
\]
Suppose now that $\Gamma\backslash G/K$ is {\em not orientable}.  Then
$H^3(\Gamma\backslash G/K, \mathbb{Z}) = \{0\}$ and therefore
$\delta(B_\sigma) = \delta(\sigma)
= 0$. Hence $B_\sigma$ is Morita equivalent to
$C(\Gamma\backslash G/K)$ and we have $(*)$ in this case.

Suppose next that $\Gamma\backslash G/K$ is {\em orientable}.
The short exact sequence of coefficient groups
\[
   1\to\mathbb{Z} \overset{i}{\to} \mathbb{R}
      \overset{e^{2\pi\sqrt{-1}}}{\longrightarrow} U(1) \to 1
\]
gives rise to a long exact sequence of cohomology groups (the change of
coefficient groups sequence)
\begin{equation}
\tag{$**$}
   \cdots \to H^2(\Gamma,\mathbb{R})
      \overset{e^{2\pi\sqrt{-1}_*}}{\longrightarrow} H^2(\Gamma, U(1))
      \overset{\delta}{\to} H^3(\Gamma,\mathbb{Z}) \overset{i_*}{\to}
      H^3(\Gamma,\mathbb{R}) \to \cdots
\end{equation}
Since $\Gamma\backslash G/K$ is oriented, we see that
$H^3(\Gamma,\mathbb{Z})\cong \mathbb{Z}$ and
$H^3(\Gamma,\mathbb{R})\cong\mathbb{R}$ are both generated by the
fundamental orientation class of $\Gamma\backslash G/K$,\ \ $
[\Gamma\backslash G/K]$, \ \ and since
\ \ $i_*[\Gamma\backslash G/K] = $ \\ $[\Gamma\backslash G/K]$, we see that
$i_*$ is
injective.
Therefore by the exactness of $(**)$ at $H^3(\Gamma, \mathbb Z)$, one has
$\delta(\sigma)=0$ for all $\sigma\in H^2(\Gamma,U(1))$, and so we see
that $B_\sigma$ is Morita equivalent to $C(\Gamma\backslash G/K)$,
and again we have $(*)$ in this case.
\end{proof}

\begin{cor} 
Let $M=K(\Gamma,1)$ be a connected locally-symmetric, compact,
3-dimensional manifold.  Let $\sigma\in H^2(\Gamma, U(1))$ be any
multiplier on $\Gamma$, then one has
\[
   K_j(C^*_r(\Gamma, \sigma)) \cong K_j(C^*_r(\Gamma)) \cong K^{j+1}(M),
      \quad j = 0,1.
\]

\end{cor}

\begin{proof}
Since $M$ is locally symmetric, it is of the form $\Gamma\backslash G/K$,
where $G$ is a connected Lie group, $K$ is a maximal compact subgroup such
that $\dim(G/K) = 3$ and $\Gamma\subset G$ is a uniform lattice in $G$.  We
need to verify that $\gamma_G=1$.  According to Thurston's list of
3-dimensional geometries or locally homogeneous spaces, one has
\begin{enumerate}
\item $G = \mathbb{R}^3 \rtimes \so(3),\ G/K = \mathbb{R}^3,\ \gamma_G=1$
since $\mathbb{R}^3$ and $\so(3)$ are amenable, and so is their semidirect
product.
\item $G=\so_0(3,1),\ G/K=\mathbb{H}^3,\ \gamma_G = 1$ by Kasparov's
theorem.
\item $G=\so_0(2,1)\rtimes\mathbb{R},\ G/K = \mathbb{H}^2\times\mathbb{R},
\gamma_G = 1$ since it's the semidirect product of $K$-amenable groups.
\item $G= \Heis,\ G/K = \Heis,\ \gamma_G = 1$ since $\Heis$ is nilpotent
and hence an amenable group.
\item $G=\Solv,\ G/K=\Solv,\ \gamma_G = 1$ since $\Solv$ is a solvable group
and hence an amenable group.
\item $G=\widetilde{\so_0(2,1)}\rtimes\mathbb{R},\
G/K=\widetilde{\so_0(2,1)}$.  Firstly, $\gamma_{\widetilde{\so_0(2,1)}}=1$
since $\widetilde{\so_0(2,1)}$ is the semidirect product of the $K$-amenable
groups $\so_0(2,1)$ and $\mathbb{Z}$.  Also $\gamma_G=1$, since its the
semidirect product of the $K$-amenable groups $\widetilde{\so_0(2,1)}$ and
$\mathbb{R}$.
\end{enumerate}

The other two locally homogeneous spaces in Thurston's list are not locally
symmetric.  We now apply Corollary 9 to deduce Corollary 10.
\end{proof}

An interesting question is whether Corollary 10 is true without the
locally symmetric assumption on $M$.  We formulate this in terms of a
conjecture.

\begin{conj*}
Let $M=K(\Gamma,1)$ be a connected, compact, 3-dimensional manifold which is
an Eilenberg-Maclane space with fundamental group $\Gamma$.  Then
for any multiplier $\sigma\in H^2(\Gamma, U(1))$ on $\Gamma$, one has
\[
   K_j(C_r^*(\Gamma,\sigma)) \cong K_j(C_r^*(\Gamma)) \cong K^{j+1}(M),
      \quad j=0,1.
\]
\end{conj*}

\begin{rems*}
Selected portions of our proof of Corollary 9 go through in the situation
described  in the conjecture.  More precisely, the proof of Corollary 9 shows
that the Dixmier-Douady invariant $\delta(\sigma) = 0$  for all $\sigma \in
H^2(\Gamma, U(1))$ for $\Gamma$ as in the conjecture.
\end{rems*}

\section{Range of the trace and the Kadison constant}

In this section, we will prove some structural theorems for the twisted
group $C^*$-algebras that are relevant to the `Martini' problems
described in the introduction. The first of these calculates the range of the
canonical trace map on $K_0$ of the twisted group $C^*$-algebras.
We use
in an essential way the results of the previous section as well as a twisted
version of the $L^2$-index theorem of Atiyah \cite{At}, which is due to
Gromov \cite{Gr2}.
This enables us to deduce information about projections in the
twisted
group $C^*$-algebras. In the case of no twisting, this follows
because the Baum-Connes conjecture is known to be true while
these results are also well
known for the case of the irrational rotation algebras. However, our
approach here is novel, and
as we will show elsewhere \cite{Ma}, enables a
generalization of most of the known results.

\subsection{Twisted Kasparov map}
Suppose that $\Gamma_g$ is a discrete, cocompact subgroup of $\so_0(2,1)$
that is, $\Gamma_g$ is the fundamental group of Riemann surface $\Sigma_g$
of genus $g>1$.  Then for any $\sigma\in H^2(\Gamma_g, U(1))$, the
\emph{twisted Kasparov isomorphism},
\begin{equation}
\tag{$*$} \mu_\sigma : K_\bullet (\Sigma_g) \to K_\bullet
(C^*_r(\Gamma_g,\sigma))
\end{equation}
is defined as follows.  Here $K_0(\Sigma_g)$ denotes the
$K$-homology group of $\Sigma_g$.  Since $\Sigma_g$ is spin, it is
$K$-oriented and by Poincar\'{e} duality, the $K$ groups
$K^j(\Sigma_g)$ are naturally isomorphic to the corresponding  $K$-homology
groups $K_j(\Sigma_g)$ for $j=0,1$.  Explicitly, let
$\mathcal{E}\to\Sigma_g$ be a
vector
bundle over $\Sigma_g$ defining  an element $[\mathcal{E}]$
in $K^0(\Sigma_g)$. Under Poincar\'{e} duality, $[\mathcal{E}]$
corresponds to the twisted Dirac
operator $\npartial^+_{\mathcal{E}} : L^2(\Sigma_g,
\mathcal{S}^+\otimes\mathcal{E})
\to L^2(\Sigma_g, \mathcal{S}^-\otimes\mathcal{E})$ where
$\mathcal{S}^{\pm}$ denote the $\frac12$
spinor bundles over $\Sigma_g$.  That is,
\begin{align*}
   PD : K^0(\Sigma_g) & \to K_0(\Sigma_g) \\
   [\mathcal{E}] & \to [\npartial_{\mathcal{E}}^+]
\end{align*}
is the Poincar\'{e} duality isomorphism. By Corollary 8 of the previous section,
there is a canonical isomorphism
$$
K_\bullet (C^*_r(\Gamma_g, \sigma)) \cong K^\bullet (\Sigma_g).
$$
Both of these maps are assembled to yield the twisted Kasparov map as in (*).

We next describe this map more explicitly.
Given $[\npartial^+_{\mathcal{E}}] \in K_0(\Sigma_g)$ as above, the lift of
this operator to $\mathbb{H} = \widetilde\Sigma_g$, the universal cover of
$\Sigma_g$,
\[
   \widetilde{\npartial_{\mathcal{E}}^+} : L^2(\mathbb{H},
      \widetilde{\mathcal{S}^+\otimes\mathcal{E}}) \to L^2(\mathbb{H},
      \widetilde{\mathcal{S}^-\otimes\mathcal{E}})
\]
is a $\Gamma_g$-invariant operator. Consider now the short
exact sequence of coefficient groups
\[
   1\to\mathbb{Z} \overset{i}{\to} \mathbb{R}
      \overset{e^{2\pi\sqrt{-1}}}{\longrightarrow} U(1) \to 1,
\]
which gives rise to a long exact sequence of cohomology groups (the change of
coefficient groups sequence)
\begin{equation}
\tag{$**$}
   \cdots \to H^2(\Gamma_g,\mathbb{Z}) \overset{i_*}{\to}
H^2(\Gamma_g,\mathbb{R})
      \overset{{e^{2\pi\sqrt{-1}}}_*}{\longrightarrow} H^2(\Gamma_g, U(1))
      \overset{\delta}{\to} 0.
\end{equation}
Therefore for any multiplier $\sigma \in H^2(\Gamma_g, U(1))$ of $\Gamma_g$,
there is a 2-form $\omega$ on $\Sigma_g$ such that
${e^{2\pi\sqrt{-1}}}_*([\omega]) = \sigma$.
Of course, the choice of $\omega$ is not unique, but this will not affect the
results that we are concerned with. Let $\widetilde \omega$ denote the lift
of $\omega$ to the universal cover $\mathbb H$. Since the hyperbolic plane
$\mathbb H$ is contractible, it follows that $\widetilde \omega = d\eta$
where $\eta$ is a 1-form on $\mathbb H$ which is not in general $\Gamma_g$
invariant. Now let $\nabla = d - i\eta$ denote a connection on the trivial
complex line bundle on $\mathbb H$. Note that the curvature of $\nabla$ is
$\nabla^2 = i\omega$. Consider now the operator
\[
   \widetilde{\npartial_{\mathcal{E}}^+}\otimes\nabla : L^2(\mathbb{H},
      \widetilde{\mathcal{S}^+\otimes\mathcal{E}}) \to L^2(\mathbb{H},
      \widetilde{\mathcal{S}^-\otimes\mathcal{E}}).
\]
It does not commute with the $\Gamma_g$ action, but it does commute
with the projective action of $\Gamma_g$ which is defined by the multiplier
$\sigma$, and by a mild generalization of the index theorem of [CM], it has a
$\Gamma_g$-$L^2$-index
\[
   \ind_{\Gamma_g} (\widetilde{\npartial^+_{\mathcal{E}}}\otimes \nabla) \in
      K_0(\mathbb{C}(\Gamma_g, \sigma) \otimes \mathcal{R}),
\]
where $\mathcal{R}$ denotes the algebra of smoothing operators.
Then observe that the {\em twisted Kasparov map} is merely
$$
\mu_\sigma([\npartial_{\mathcal{E}^+}]) = j_*(\ind_{\Gamma_g}
(\widetilde{\npartial^+_{\mathcal{E}}}\otimes \nabla) ) \in
K_0({C}^*(\Gamma_g, \sigma)),
$$
where $j: {\mathbb C}(\Gamma_g, \sigma)\otimes \mathcal R \to
{C}^*_r(\Gamma_g, \sigma)\otimes \mathcal K$
is the natural inclusion map, and
$$j_* : K_0 ( {\mathbb C}(\Gamma_g,
\sigma)\otimes \mathcal R) \to K_0 ({C}^*_r(\Gamma_g, \sigma))$$ is the
induced map on $K_0$.

The canonical trace on ${C}^*_r(\Gamma_g, \sigma))$ induces a linear map
$$
[\tr] : K_0 ({C}^*_r(\Gamma_g, \sigma)) \to \mathbb R
$$
which is called the {\em trace map} in $K$-theory.
Explicitly, first $\tr$ extends to matrices with entries in
${C}^*(\Gamma_g, \sigma)$ as (with Trace denoting matrix trace):
\[
   \tr(f\otimes r) = {\mbox{Trace}}(r) \tr(f).
\]

Then the extension of $\tr$ to $K_0$ is given by
$[\tr]([e]-[f]) = \tr(e) - \tr(f)$, where $e,f$ are idempotent matrices with
entries in ${C}^*(\Gamma_g, \sigma))$.

\subsection{The isomorphism classes of algebras ${C}^*(\Gamma_g, \sigma)$}
Let $\sigma \in Z^2(\Gamma_g, U(1))$ be a multiplier on $\Gamma_g$. If
$\sigma' \in Z^2(\Gamma_g, U(1))$ is another multiplier on $\Gamma_g$ such
that $[\sigma] = [\sigma'] \in H^2(\Gamma_g, U(1))$, then it can be easily
shown that ${C}^*(\Gamma_g, \sigma) \cong {C}^*(\Gamma_g, \sigma')$. That is,
the isomorphism classes of the $C^*$-algebras ${C}^*(\Gamma_g, \sigma)$ are
naturally parametrized by $H^2(\Gamma_g, U(1))$. But $H^2(\Gamma_g, U(1)) \cong
H^2(\Sigma_g, U(1)) \cong U(1)$ and the isomorphism is given explicitly by
$[\sigma] \to
<[\sigma], [\Sigma_g]>$, where $[\sigma]$ is now viewed as a \v{C}ech 2-cocycle
on $\Sigma_g$ with coefficients in $U(1)$, and $[\Sigma_g]$ denotes the
fundamental
class of the genus $g$ Riemann surface. We summarize this below.

\begin{lemma} The isomorphism classes of twisted group $C^*$-algebras
${C}^*(\Gamma_g, \sigma)$ are
naturally parametrized by $U(1) \cong \mathbb R/\mathbb Z \cong (0,1]$. The
classification map  is given explicitly by
$$[\sigma] \to
<[\sigma], [\Sigma_g]>,$$
 where $[\sigma]$ is now viewed as a \v{C}ech 2-cocycle
on $\Sigma_g$ with coefficients in $U(1)$, and $[\Sigma_g]$ denotes the
fundamental
class of the genus $g$ Riemann surface.
\end{lemma}

\subsection{Range of the trace map on $K_0$}
We can now state the first major theorem of this section.

\begin{thm}
The range of the trace map is
$$
[\tr] (K_0 ({C}^*_r(\Gamma_g, \sigma)) ) = \mathbb Z \theta + \mathbb Z,
$$
where $2\pi\theta = <\sigma, [\Sigma_g]>\ \in (0,1]$ is the result of
pairing the
multiplier $\sigma$ on $\Gamma_g$ with the fundamental class of $\Sigma_g$.
\end{thm}

\begin{proof}

We first observe that by the results of the previous section the
twisted Kasparov map is an isomorphism. Therefore to compute the range of the
trace map on $K_0$, it suffices to compute the range of the trace map
on elements of the form
$$\mu_\sigma([\npartial_{\mathcal{E}^0 }^+] -
[\npartial_{\mathcal{E}^1 }^+])$$
 for any element
 $$[\npartial_{\mathcal{E}^0}^+] -
[\npartial_{\mathcal{E}^1 }^+] \in K_0(\Sigma_g).$$

By the twisted analogue of the $L^2$ index theorem of Atiyah \cite{At} and
Singer \cite{Si} for elliptic operators on a covering space that are invariant
under the projective action of the fundamental group defined by $\sigma$,
and which is due to Gromov \cite{Gr2} (see also \cite{Ma}
for a detailed proof of a further generalization),
one has
\begin{equation}
\tag{$*$}
   [\tau](\ind_{\Gamma_g} (\widetilde{\npartial_{\mathcal{E}}^+}\otimes
\nabla)) =
       \frac{1}{2\pi}\langle \hat{A}(\Sigma_g) \ch(\mathcal{E}) e^{[\omega]}
      , [\Sigma_g]\rangle.
\end{equation}
We next simplify the right hand side of $(*)$ using
\begin{align*}
   \hat{A}(\Sigma_g) &= 1  \\
   \ch(\mathcal{E}) &= \rank \mathcal{E} + c_1(\mathcal{E}) \\
   e^{[\omega]} &= 1 + {[\omega]} .
\end{align*}
Therefore one has
\[
   [\tau](\ind_{\Gamma_g} (\widetilde{\npartial_{\mathcal{E}}^+}\otimes
\nabla)) =
      {\rank \mathcal E}\frac{\langle [\omega], [\Sigma_g] \rangle}{2\pi}
      + \frac{\langle c_1(\mathcal{E}), [\Sigma_g] \rangle}{2\pi},
\]
and we see that
$$
[\tau](\mu_\sigma([\npartial_{\mathcal{E}^0 }^+] -
[\npartial_{\mathcal{E}^1 }^+])) = ({\rank {\mathcal E}^0 - {\rank
{\mathcal E}^1})}
\frac{\langle [\omega], [\Sigma_g] \rangle}{2\pi}
      + \frac{\langle c_1(\mathcal{E}^0) - c_1(\mathcal{E}^1), [\Sigma_g]
\rangle}{2\pi}.
$$
It follows that the range of the trace map on $K_0$ is
$\mathbb Z \frac{\langle [\omega], [\Sigma_g] \rangle}{2\pi}  + \mathbb Z =
\mathbb Z \theta + \mathbb Z$,
because $$\frac{\langle [\omega], [\Sigma_g] \rangle}{2\pi} -
\theta\in \mathbb Z.$$
\end{proof}

We will now discuss some applications of this result. We begin by studying
projections in the twisted group $C^*$-algebra, which is a problem of
 independent
interest.

\begin{prop}
Let $\sigma \in H^2(\Sigma_g, \mathbb R/ \mathbb Z)$ be a multiplier on
$\Gamma_g$,  and
$2\pi\theta = <\sigma, [\Sigma_g]> \in (0,1]$ be the result of pairing
$\sigma$ with the fundamental class of $\Sigma_g$. If $\theta = p/q$ is
rational, then
there are only $q-1$ unitary equivalence classes of
projections, other than $0$ and $1$,
in the reduced twisted group $C^*$-algebra ${C}^*_r(\Gamma_g, \sigma)$.
\end{prop}

\begin{proof} By assumption, $\theta = p/q$.
Let $P$ be a projection in ${C}^*_r(\Gamma_g, \sigma)$. Then
$1-P$ is also a projection in ${C}^*_r(\Gamma_g, \sigma)$ and one has
$$
1= \tr(1) = \tr(P) + \tr(1-P).
$$
Each term in the above equation is non-negative. By the previous
theorem,
it follows that $\tr(P) \in \{0, 1/q, 2/q, \ldots 1\}$. By faithfulness
and normality  of the
trace $\tr$, it follows that there are only $q-1$
unitary equivalence classes of projections, other than those of
$0$ and $1$
in ${C}^*_r(\Gamma_g, \sigma)$.
\end{proof}

Our second application will involve the Kadison constant of a twisted group
$C^*$-algebra, which we will now recall.
The {\em Kadison constant} of ${C}^*_r(\Gamma_g, \sigma)$ is defined by:
$$
C_\sigma(\Gamma_g) = \inf\{ \tr(P) : P \ \ {\mbox{is a
non-zero projection in}} \ \
{C}^*_r(\Gamma_g, \sigma) \otimes \mathcal K\}.
$$
Recall from earlier sections the following Hamiltonians:
$$
H_\eta = (d-i\eta)^*(d-i\eta) = \nabla^*\nabla,
$$
and
$$
H_{\eta, V} = H_\eta + V,
$$
where $V$ is any $\Gamma_g$-invariant potential on $\mathbb H$. The
operators
$H_\eta$ and $H_{\eta, V}$ are invariant under the projective $(\Gamma_g,
\sigma)$-action.

\begin{prop}
Let $\sigma \in H^2(\Sigma_g, \mathbb R/ \mathbb Z)$ be a multiplier on
$\Gamma_g$,  and
$2\pi\theta = <\sigma, [\Sigma_g]> \in (0,1]$ be the result of pairing
$\sigma$ with the fundamental class of $\Sigma_g$. If $\theta = p/q$ is
rational, then
the spectrum of any associated Hamiltonian $H_{\eta, V}$ has a band
structure, in the
sense that the intersection of the resolvent set with any compact interval
in $\mathbb R$
has only a finite number of components. In particular, the intersection of
$\sigma(H_{\eta, V})$ with any compact interval in $\mathbb R$ is never a
Cantor set.
\end{prop}

\begin{proof}
By the previous proposition, it follows that one has the estimate
$C_\sigma(\Gamma_g)
\ge 1/q >0$. Then one applies the main result in Br\"uning-Sunada
\cite{BrSu} to
deduce the proposition.
\end{proof}

This leaves open the question of whether there are Hamiltonians
with Cantor spectrum
when $\theta$ is irrational. In the Euclidean case, this is usually known
as the {\em Ten Martini Problem}, and is to date, not completely solved,
though much progress has been made (cf. \cite{Sh}). We pose a
generalization
of this problem to the hyperbolic case (which also includes the
Euclidean case):

\begin{conj*}[The Ten Dry Martini Problem]
Let $\sigma \in H^2(\Sigma_g, \mathbb R/ \mathbb Z)$ be a multiplier on
$\Gamma_g$,  and
$2\pi\theta = <\sigma, [\Sigma_g]> \in (0,1]$ be the result of pairing
$\sigma$ with the fundamental class of $\Sigma_g$. If $\theta$ is
irrational, then
there is an associated Hamiltonian $H_{\eta, V}$ with a Cantor set
type spectrum, in the sense that the intersection of
$\sigma(H_{\eta, V})$ with some compact interval in $\mathbb R$ is a Cantor set.
\end{conj*}

We will next apply the range of the trace Theorem 12 to deduce results about
the {\em discrete Hamiltonian} $H_\tau$, as in Section 9.

\begin{prop}
Let $\sigma \in H^2(\Sigma_g, \mathbb R/ \mathbb Z)$ be a multiplier on
$\Gamma_g$,  and
$2\pi\theta = <\sigma, [\Sigma_g]> \in (0,1]$ be the result of pairing
$\sigma$ with the fundamental class of $\Sigma_g$. If $\theta = p/q$ is
rational, then
the spectrum of the associated discrete Hamiltonian $H_{\tau}$ has a band
structure, in the
sense that the intersection of the resolvent set with $\mathbb R$
has only a finite number of components. In particular, the intersection of
$\sigma(H_{\tau})$ with any compact interval in $\mathbb R$ is never a
Cantor set.
\end{prop}

\begin{proof} From the estimate
$C_\sigma(\Gamma_g)
\ge 1/q >0$ the main result in \cite{Sun} implies
the proposition.
\end{proof}
This leads us to our next conjecture.

\begin{conj*}[The Discrete Ten Dry Martini Problem]
Let $\sigma \in H^2(\Sigma_g, \mathbb R/ \mathbb Z)$ be a multiplier on
$\Gamma_g$,  and
$2\pi\theta = <\sigma, [\Sigma_g]> \in (0,1]$ be the result of pairing
$\sigma$ with the fundamental class of $\Sigma_g$. If $\theta$ is
irrational, then
the associated Hamiltonian $H_{\tau}$ has Cantor spectrum.
\end{conj*}

\subsection{On the classification of twisted group $C^*$-algebras}
We will now use the range of the trace Theorem 12, to give a complete
classification,
up to isomorphism, of the twisted group $C^*$-algebras ${C}^*(\Gamma,
\sigma)$. A similar complete classification, up to Morita equivalence,
is contained in \cite{Ma}.

\begin{prop}[Isomorphism classification of twisted group $C^*$-algebras]
Let $\sigma, \sigma' \in $ \\ $H^2(\Sigma_g, \mathbb R/ \mathbb Z)$ be
multipliers on
$\Gamma_g$, and $2\pi\theta = <\sigma, [\Sigma_g]> \in (0,1]$,
$2\pi\theta' = <\sigma', [\Sigma_g]> \in (0,1]$
be the result of pairing
$\sigma$, $\sigma'$ with the fundamental class of $\Sigma_g$.
Then ${C}^*(\Gamma_g, \sigma) \cong {C}^*(\Gamma_g, \sigma')$ if and only if
$\theta' \in \{\theta, 1-\theta\}$.
\end{prop}

\begin{proof}
 Let $\tr$ and $\tr'$ denote the canonical traces on ${C}^*(\Gamma_g, \sigma)$
and ${C}^*(\Gamma_g, \sigma')$ respectively. Let
$$
\phi : {C}^*(\Gamma_g, \sigma) \cong {C}^*(\Gamma_g, \sigma')
$$
be an isomorphism, and let
$$
\phi_* : K_0({C}^*(\Gamma_g, \sigma)) \cong K_0({C}^*(\Gamma_g, \sigma'))
$$
denote the induced map on $K_0$.
By Theorem 12, the range of the trace map on $K_0$ is
$$
[\tr] (K_0 ({C}^*(\Gamma_g, \sigma)) ) = \mathbb Z \theta + \mathbb Z
$$
and
$$
[\tr'] (K_0 ({C}^*(\Gamma_g, \sigma')) ) = \mathbb Z \theta' + \mathbb Z.
$$
So there are elements $[P] \in K_0({C}^*(\Gamma_g, \sigma))$
and $[P'] \in K_0({C}^*(\Gamma_g, \sigma') )$ such that
$[\tr]([P]) = \theta$ and $[\tr']([P']) = \theta'$.
Clearly one has $\tr\circ \phi = \tr'$,
which induces the identity $[\tr]\circ \phi_* = [\tr']$ in
$K_0 ({C}^*(\Gamma_g, \sigma'))$.
In Section 1, we have proved that $K_0 ({C}^*(\Gamma_g, \sigma)) \cong
\mathbb Z^2 \cong K_0 ({C}^*(\Gamma_g, \sigma'))$. In the basis above,
one has
$$
\phi_* : \mathbb Z[P] \oplus \mathbb Z \cong K_0({C}^*(\Gamma_g, \sigma))
\to K_0({C}^*(\Gamma_g, \sigma')) \cong \mathbb Z[P'] \oplus \mathbb Z.
$$
Since $\phi_*[1] = [1]$ and $\phi_* \in {\mathbf{GL}}(2, \mathbb Z)$, one
sees that
there is an
integer $n$ such that
$$
\phi_* = \left(
\begin{array}{cc}
1 & n \\[2mm]
0 & \pm 1 \end{array}
\right).
$$
Assembling these results, one has
$ \theta = [\tr]([P]) = [\tr](\phi_*[P]) = [\tr'] (n[1] \pm [P'])
= n \pm\theta'$.
Since $\theta, \theta' \in (0,1]$, one deduces that $\theta' \in \{\theta,
1-\theta\}$.

Let $\psi : \Sigma_g \to \Sigma_g$ be an orientation reversing diffeomorphism.
We can assume without loss of generality that $\psi$ has a fixed point
$x_0\in \Sigma_g$.
This is because there is an orientation preserving diffeomorphism $\eta$
of $\Sigma_g$ whose value at the
point $\psi(x_0)$ is equal to $x_0$; in fact $\eta$ can be chosen to be
isotopic to
the identity (cf. exercise A2, chapter 1, \cite{Helg}). The composition
$\eta\circ \psi$
is then an orientation reversing diffeomorphism of $\Sigma_g$ with fixed
point $x_0$.
Then $\psi$ induces an automorphism $\psi_*:\Gamma_g \to \Gamma_g$ of the
fundamental group
$\pi_1(\Sigma_g, x_0) \cong \Gamma_g$. We first evaluate
$<\psi^*\sigma, [\Sigma_g]> = <\sigma, \psi_*[\Sigma_g]>
={\overline{<\sigma, [\Sigma_g]>}}=
<\bar\sigma, [\Sigma_g]>$,
since $\psi$ is orientation
reversing. By Lemma 13  we see that $\psi^*\sigma =
\bar\sigma\in H^2(\Gamma_g, U(1))$.
Therefore the automorphism $\psi_*$ of $\Gamma_g$ induces an isomorphism of
twisted
group $C^*$-algebras
$$
{C}^*(\Gamma_g, \sigma) \cong {C}^*(\Gamma_g, \psi^*\sigma)\cong
{C}^*(\Gamma_g, \bar\sigma).
$$
Therefore if
$\theta' \in \{\theta, 1-\theta\}$, one has
${C}^*(\Gamma_g, \sigma) \cong {C}^*(\Gamma_g, \sigma')$, completing the
proof of
the proposition.
\end{proof}

\subsection{Twisted ICC group von Neumann algebras and type II$_1$ factors}
Recall that an ICC group $\Gamma$ is one in which every non-trivial
conjugacy class  is infinite.
There are many examples of ICC groups, such as free groups, fundamental groups
of compact surfaces etc. It is well known that the group  von Neumann algebras
of these groups are type II$_1$ factors \cite{Tak}. We will now prove that a
similar result holds for the twisted group von Neumann algebras
(this result probably exists in the literature but
for completeness we reproduce a proof). We briefly
recall some definitions. Let $W^*(\Gamma, \sigma)$ denote the twisted group
von Neumann algebra, where $\sigma$ is a multiplier on $\Gamma$, which is
by definition the {\em weak} closure of
${C}^*(\Gamma, \sigma)$, or equivalently, the weak closure of
the algebraic group algebra
$\mathbb{C}(\Gamma, \sigma)$ in the
$\sigma$-regular representation on
 $\ell^2(\Gamma)$. Let
 ${\mbox{Proj}}(W^*(\Gamma,\sigma))$
denote the set of all projections in $W^*(\Gamma, \sigma)$. Then one has.

\begin{prop}
Let $\Gamma$ be an ICC group, and
$\sigma \in H^2(\Sigma, \mathbb R/ \mathbb Z)$ be a multiplier on $\Gamma$.
Then $W^*(\Gamma, \sigma)$ is a II$_1$ factor. In particular,
$\tr( {\mbox{Proj}}(W^*(\Gamma, \sigma))) = [0,1]$.
\end{prop}

\begin{proof}
By
the commutant theorem for the regular $\sigma$-representation
we see that the commutant of $W^*(\Gamma, \sigma)$ is identified with
$W^*(\Gamma, \bar\sigma)$. We need to compute the centre $Z(\Gamma, \sigma)$
of $W^*(\Gamma, \sigma)$, which is equal to the intersection
$Z(\Gamma, \sigma) = W^*(\Gamma, \sigma) \cap W^*(\Gamma, \bar\sigma)$.
Let $T : \Gamma \to B(\ell^2(\Gamma))$ denote the left projective $(\Gamma,
\sigma)$-action.
Regard $x\in W^*(\Gamma, \bar\sigma)$ as $x = \sum_{\gamma\in
\Gamma} x(\gamma) T(\gamma)$.
Since $W^*(\Gamma, \sigma)$ is  the {\em weak} closure of
$\mathbb{C}(\Gamma, \sigma)$, it follows that $ (x(\gamma))_{\gamma\in\Gamma}
 \in \ell^2(\Gamma)$.
Now $x\in Z(\Gamma, \sigma)$ if and only if $x$ commutes with
$T(\gamma'), \ \gamma'\in \Gamma$. But
\begin{align*}
T(\gamma') x  & = \sum_{\gamma\in \Gamma} x(\gamma) \sigma(\gamma',\gamma)
T(\gamma'\gamma)\\\
              & = \sum_{\gamma\in \Gamma} x(\gamma^{\prime^{-1}}\gamma)
\sigma(\gamma',\gamma^{\prime^{-1}}\gamma) T(\gamma),
\end{align*}
and
\begin{align*}
 x T(\gamma') & = \sum_{\gamma\in \Gamma} x(\gamma) \sigma(\gamma,\gamma')
T(\gamma\gamma')\\\
              & = \sum_{\gamma\in \Gamma} x(\gamma\gamma^{\prime^{-1}})
\sigma(\gamma\gamma^{\prime^{-1}},\gamma') T(\gamma).
\end{align*}
Therefore we see that
$x(\gamma^{\prime^{-1}}\gamma) \sigma(\gamma',\gamma^{\prime^{-1}}\gamma)=
x(\gamma\gamma^{\prime^{-1}}) \sigma(\gamma\gamma^{\prime^{-1}},\gamma')$
for all $\gamma'\in
\Gamma$. So
$|x(\gamma^{\prime^{-1}}\gamma\gamma')| = |x(\gamma)|$ for all $\gamma'\in
\Gamma$. That is,
$|x(\cdot)|$ is constant on each conjugacy
class. Now since $x\in \ell^2(\Gamma)$, it follows that $x$ vanishes
on each infinite
conjugacy class. Since $\Gamma$ is an ICC group, it follows that $x(\gamma) = 0$
for all $\gamma\ne 1$, that is $Z(\Gamma, \sigma)$ is 1-dimensional and
$W^*(\Gamma, \sigma)$ is a II$_1$ factor.
\end{proof}

\section{The topological index and the index theorem}

This section  identifies the  Hall conductivity
$\tau_c(P,P,P) = \tau(P\, dP\, dP)$ with a topological invariant,
generalising the work of \cite{Xia}.
Suppose that $\Gamma_g$ is a discrete, cocompact subgroup of $\so_0(2,1)$.
That is, $\Gamma_g$ is the fundamental group of a Riemann surface $\Sigma_g$
of genus $g>1$.  Then for any $\sigma\in H^2(\Gamma_g, U(1))$, the
\emph{twisted Kasparov isomorphism},
\[
   \mu_\sigma : K_\bullet (\Sigma_g) \to K_\bullet (C^*_r(\Gamma_g,\sigma))
\]
is defined as in the previous
section.
We note in the following section  (using a
result of \cite{Ji})
that given any
projection $P$ in $C^*_r(\Gamma,\sigma)$ there is
both a projection $\tilde P$ in the same $K_0$ class but lying in
a dense subalgebra, stable under the holomorphic functional calculus,
and a Fredholm module for this dense subalgebra, which may be
be paired with $\tilde P$
to obtain an analytic index. On the other hand, by the results of the current
section,
given any such projection $P$ there
is a topological index that we can associate to it. The main result we prove
here is that the (analytic index) = (topological index).

The first step in the proof is to show that
given an additive
 group cocycle $c\in Z^2(\Gamma_g)$ we may define canonical pairings
with $K_0(\Sigma_g)$ and $K_0(C^*_r(\Gamma_g, \sigma))$ which
are related by the twisted Kasparov isomorphism, by
generalizing some of the results of Connes and
Connes-Moscovici to the twisted case.
The group 2-cocycle $c$ may be
regarded as a skew symmetrised function on
$\Gamma_g\times\Gamma_g\times\Gamma_g$, so that
 we can modify a standard
construction in \cite{CM} to obtain a cyclic
2-cocycle $\tau_c$ on $\mathbb{C}(\Gamma_g, \sigma)\otimes \mathcal{R}$
by defining:
\[
   \tau_c(f^0\otimes r^0, f^1\otimes r^1, f^2\otimes r^2) =
{\mbox{Tr}}(r^0r^1r^2)
      \sum_{g_0g_1g_2=1} f^0(g_0) f^1(g_1) f^2(g_2) c(1, g_1, g_1g_2)
\sigma(g_1, g_2).
\]
Note that
 $\tau_c$ extends to $\mathbb{C}(\Gamma_g, \sigma) \otimes \mathcal{L}^2$,
(where
$\mathcal{L}^2$ denotes Hilbert-Schmidt operators) and by the pairing
theory of \cite{Co} one gets an additive map
\[
   [\tau_c] : K_0(\mathbb{C}(\Gamma_g, \sigma) \otimes
\mathcal{R})\to\mathbb{R}.
\]
Explicitly, $[\tau_c]([e]-[f]) = \widetilde\tau_c(e,\cdots,e) -
\widetilde\tau_c(f,\cdots,f)$, where $e,f$ are idempotent matrices with
entries in $(\mathbb{C}(\Gamma_g, \sigma)\otimes \mathcal{R})^\sim = $ unital
algebra obtained by adding the identity to
$\mathbb{C}(\Gamma_g, \sigma)\otimes\mathcal{R}$ and $\widetilde\tau_c$
denotes the
canonical extension of $\tau_c$ to $(\mathbb{C}(\Gamma_g, \sigma)\otimes
\mathcal{R})^\sim$. Let \
$\widetilde{\npartial_{\mathcal{E}}^+}\otimes \nabla$ \  be
the Dirac operator defined in the previous section, which is
invariant under the projective action of the fundamental group
defined by $\sigma$.
By definition, the $(c,\Gamma_g, \sigma)$-index of
\ $\widetilde{\npartial_{\mathcal{E}}^+}\otimes \nabla$\ is
\[
   [\tau_c](\ind_{\Gamma_g} (\widetilde{\npartial_{\mathcal{E}}^+}\otimes
\nabla))\in\mathbb{R}.
\]
It only depends on the cohomology class $[c]\in H^2(\Gamma_g)$, and it is
linear with respect to $[c]$.  We assemble this to give the following theorem.

\begin{thm}
Given $[c] \in H^2(\Gamma_g)$ and $\sigma \in H^2(\Gamma_g, U(1))$ a
multiplier on $\Gamma_g$, there is a canonical additive map
\[
   \langle [c],\ \rangle : K_0(\Sigma_g)\to\mathbb{R},
\]
which is defined as
\[
   \langle [c], [{\npartial_{\mathcal{E}}^+}]\rangle = [\tau_c]
      (\ind_{\Gamma_g} (\widetilde{\npartial_{\mathcal{E}}^+}\otimes
\nabla)) \in\mathbb{R}.
\]
Moreover, it is linear with respect to $[c]$.
\end{thm}

By a generalization of the Connes-Moscovici higher index theorem
\cite{CM} to the twisted case of
elliptic operators on a covering space that are invariant
under the projective action of the fundamental group defined by $\sigma$,
(see \cite{Ma} for a detailed proof),
 one has
\begin{equation}
\tag{$*$}
   [\tau_c](\ind_{\Gamma_g} (\widetilde{\npartial_{\mathcal{E}}^+}\otimes
\nabla)) =
       \frac{1}{2\pi}\langle \hat{A}(\Sigma_g) \ch(\mathcal{E}) e^{[\omega]}
      \psi^*(c), [\Sigma_g]\rangle ,
\end{equation}
where $\psi: \Sigma_g \to \Sigma_g$ is the classifying map of the universal
cover (which in this case is the identity map) and $[c]$ is considered as a
degree 2 cohomology class on $\Sigma_g$.  We next simplify the right hand
side of $(*)$ using the fact that $\hat{A}(\Sigma_g) = 1$ and that
\begin{align*}
   \ch(\mathcal{E}) &= \rank \mathcal{E} + c_1(\mathcal{E}), \\
   \psi^*(c) &= c,\\
   e^{[\omega]} &= 1 + {[\omega]} .
\end{align*}
We obtain
\[
   [\tau_c](\ind_{\Gamma_g} (\widetilde{\npartial_{\mathcal{E}}^+}\otimes
\nabla)) =
      \frac{\rank \mathcal E}{2\pi}\langle [c], [\Sigma_g] \rangle .
\]

\begin{cor}
Let $c,\ [c]\in H^2(\Gamma_g)$, be the area cocycle.  Then one has
\[
   \langle [c], [{\npartial_{\mathcal{E}}^+}]\rangle = 2(g-1)\rank\mathcal{E}
      \in\mathbb{Z}.
\]
\end{cor}

\begin{proof*}
When $c,\ [c]\in H^2(\Gamma_g)$, is the area 2-cocycle, one has
\[
   \langle [c], [\Sigma_g] \rangle = -2\pi\chi(\Sigma_g) = 4\pi(g-1).
\]
\end{proof*}

\begin{rems}
These theorems have been generalised in \cite{Ma}.
They agree with Xia's result \cite{Xia}, although our methods are different.
\end{rems}

We next describe the canonical pairing of $K_0(C^*_r(\Gamma_g, \sigma))$,
given $[c]\in H^2(\Gamma_g)$. Since $\Sigma_g$ is negatively curved, we
know from \cite{Ji} that
\[
   A_{\sigma, g} = \left\{ f : \Gamma_g \to \mathbb{C} \mid
      \sum_{\gamma\in\Gamma_g} |f(\gamma)|^2 (1+l(\gamma))^k < \infty
      \mbox{ for all } k\ge 0\right\},
\]
where $l:\Gamma_g \to \mathbb{R}^+$ denotes the length function, is a dense
and spectral invariant subalgebra of $C_r^*(\Gamma_g, \sigma)$.  In
particular it is closed under the smooth functional calculus, and is
known as the algebra of rapidly decreasing $L^2$ functions on $\Gamma_g$.
By a theorem of \cite{Bost}, the inclusion map $A_{\sigma, g}\subset
C^*_r(\Gamma_g, \sigma)$ induces an isomorphism
\[
   K_j(A_{\sigma, g}) \cong K_j(C_r^*(\Gamma_g, \sigma)),\quad j=0,1.
\]
As  $\Sigma_g$ is a negatively curved manifold,
we know (by \cite{Mos} and \cite{Gro}) that degree 2 cohomology classes in
$H^2(\Gamma_g)$ have
\emph{bounded} representatives i.e.\  bounded 2-cocycles on $\Gamma_g$.
Let $c$ be a bounded 2-cocycle on $\Gamma_g$.  Then it defines a cyclic
2-cocycle $\tau_c$ on the twisted group algebra
${\mathbb C}(\Gamma_g, \sigma)$, by a slight modification of
the standard formula \cite{CM}, (\cite{Ma} for the general case)
\[
   \tau_c(f^0, f^1, f^2) = \sum_{g_0g_1g_2=1} f^0(g_0) f^1(g_1) f^2(g_2)
      c(1, g_1, g_1g_2)\sigma(g_1, g_2).
\]
Here $c$ is assumed to be skew-symmetrized. Since the only difference
with the expression obtained in \cite{CM} is $\sigma(g_1, g_2)$, and since
$|\sigma(g_1, g_2)| = 1$, we can use Lemma 6.4, part (ii) in \cite{CM}
and the assumption that $c$ is bounded, to obtain the necessary
estimates which show that in fact $\tau_c$ extends continuously to the bigger
algebra $A_{\sigma, g}$. This induces an additive map in $K$-theory as before:
\begin{gather*}
   [\tau_c] : K_0(A_{\sigma, g})\to\mathbb{R} \\
   [\tau_c]([e] - [f]) = \widetilde\tau_c(e,\cdots, e)
      -\widetilde\tau_c(f,\cdots,f),
\end{gather*}
where $e,f$ are idempotent matrices with entries in $(A_{\sigma, g})^\sim$ (the
unital algebra associated to $A_{\sigma, g}$) and $\widetilde\tau_c$ is
the canonical extension of $\tau_c$ to $(A_{\sigma, g})^\sim$.
Observe that the twisted Kasparov map is merely
$$
\mu_\sigma([\npartial_{\mathcal{E}^+}]) = j_*(\ind_{\Gamma_g}
(\widetilde{\npartial^+_{\mathcal{E}}}\otimes \nabla) ) \in
K_0({C}^*(\Gamma_g, \sigma)).
$$
Here $j: {\mathbb C}(\Gamma_g, \sigma)\otimes \mathcal R \to
{C}^*(\Gamma_g, \sigma)\otimes \mathcal K$
is the natural inclusion map, and $j_* : K_0 ( {\mathbb C}(\Gamma_g,
\sigma)\otimes \mathcal R) \to K_0 ({C}^*(\Gamma_g, \sigma))$ is the
induced map in $K$-theory. Therefore one has the equality
\[
   \langle [c], \mu_\sigma^{-1}[P]\rangle = \langle [\tau_c], [P] \rangle
\]
for any $[P]\in K_0(A_{\sigma, g}) \cong K_0(C_r^* (\Gamma_g, \sigma))$.
Using the previous corollary, one has

\begin{cor}
Let $c,\ [c]\in H^2(\Gamma_g)$, be the area 2-cocycle.  Then $c$ is known to
be a bounded 2-cocycle, and one has
\[
   \langle [\tau_c], [P] \rangle = 2(g-1)(\rank\mathcal{E}^0 -
      \rank\mathcal{E}^1)\in\mathbb{Z},
\]
where $[P]\in K_0(A_{\sigma, g}) \cong K_0(C_r^*(\Gamma_g, \sigma))$, and where
\[
   \mu_\sigma^{-1}[P] = [{\npartial_{\mathcal{E}^0}^+}] -
      [{\npartial_{{\mathcal E}^1}^+}] \in K_0(\Sigma_g).
\]
\end{cor}

\begin{rems}
This generalizes the main result of Xia, \cite{Xia}.

We will next prove the existence of a canonical element in $
KK(C_r^*(\Gamma_g, \sigma), \mathbb{C})$, which we call the twisted
Mishchenko element.

\begin{thm}[The twisted Mishchenko element]
There exists a
unique element $[m_\sigma] \in KK(C_r^*(\Gamma_g, \sigma), \mathbb{C})$, called
the {\em twisted Mishchenko element}, such that
\begin{equation}
\tag{$*$} [1] \otimes_{C_r^*(\Gamma_g, \sigma)} [m_\sigma] = 2(g-1),
\end{equation}
where $[1] \in K_0(C_r^*(\Gamma_g, \sigma))$ denotes the module generated
by $C_r^*(\Gamma_g, \sigma)$.
\end{thm}

\begin{proof}
By the well definedness of the Kasparov intersection product \cite{Kas2},
the equation $(*)$ above defines the element $[m_\sigma]$ uniquely. In the
next section we construct a 2-summable Fredholm module $(F,\mathcal{H})$,
which defines an element $[(F,\mathcal{H})] \in KK(C_r^*(\Gamma_g, \sigma),
\mathbb{C})$, and
whose Chern character is the cyclic area 2-cocycle $[\tau_c]$, (cf.
\cite{Co2}) defined by
the area 2-cocycle $c$ on the discrete group $\Gamma$. We compute that
$$
[1] \otimes_{C_r^*(\Gamma_g, \sigma)} [(F,\mathcal{H})] = {\mbox{index}}(F)
= \tau_c(1,1,1) = 2(g-1).
$$
By uniqueness (proved above), we see that  $[m_\sigma] =  [(F,\mathcal{H})]$,
which establishes existence.
\end{proof}

This completes the proof
of Theorem 4 and Corollary 6 because we regard $\Index(PFP)$ as the result
of pairing an element of the K-homology of $\Sigma_g$ (defined by the twisted
Mishchenko element)
 with an element of
$K_0({\ck}^\Gamma )\cong K_0(C^*(\Gamma_g, \sigma))$.
This enables us to demonstrate the relationship between Corollary
12 and the discrete model of the hyperbolic Hall effect.
\end{rems}

\section{A discrete Fredholm module and the analytic index}

We have observed following Sunada that $H_\tau$ is an operator in
the twisted algebraic group algebra ${\mathbb C}(\Gamma,\sigma)$,
which is a subalgebra of $A_{\sigma,g}$.
We remark that a spectral projection
into a gap in the spectrum of $H_\tau$ is given by the smooth functional
calculus
applied to $H_\tau$. It follows from \cite{Ji} that such spectral projections
lie in $A_{\sigma,g}$.
Connes constructs a Fredholm module for
${\mathbb C}\Gamma$
 which can be adapted
to the case of
 ${\mathbb C}(\Gamma,\sigma)$.
In his construction the Hilbert space is the $\ell^2$ sections of
 the restriction of the spinor bundle to the orbit $\Gamma.u$.
This space is isomorphic to ${\mathcal H}=\ell^2(\Gamma)\oplus\ell^2(\Gamma)$
under the map $\iota\oplus\iota$.
The grading is the obvious one given by the $2\times 2$ matrix $\varepsilon$.
We may define the operator $F$
as in Section 7  to be multiplication by the matrix function
$$
\left(\begin{array}{cc}
0 & \varphi^{*} \\ \varphi & 0 \end{array}\right),
$$
where we restrict $\varphi$ to the orbit $\Gamma.u$.

Connes \cite{Co2}  shows that the module of the previous paragraph
is 2-summable
for ${\mathbb C}\Gamma$. We show below using the same argument
as in \cite{Co2} that if $\lambda$ denotes the left regular
$\sigma$ representation of $C^*(\Gamma,\sigma)$ then $[F,\lambda(\gamma)]$
is Hilbert-Schmidt.   So $({\mathcal H}, F)$
 is also
a 2-summable module for ${\mathbb C} (\Gamma,\sigma)$.
We may also exploit \cite{Co2} to determine explicitly the
character of this Fredholm module for our case.
We now summarise some of the pertinent details.

First, we are using the usual trace $\tr$ on the bounded operators
on ${\mathcal H}$. Second, our module
is the $\ell^2$ sections of the
 restriction of the spinor bundle to the orbit.
>From this point of view
$F$ corresponds to Clifford multiplication of a unit tangent
vector to a geodesic connecting
a given vertex of the graph to a point $x_0\notin\Gamma.u$.
We use the same notation $\varphi(\gamma.u)$  for
this unit tangent vector, regarding $\varphi$ as a  function from
$\Gamma.u$ to $T(\hyp )$, the tangent space of $\hyp$, as no confusion
will arise.

Next, note that for $f\in \mathcal H$,
 $$ [F,\lambda(\gamma)]f(\gamma') =
 (\varphi(\gamma.u)'-\varphi(\gamma^{-1}\gamma'.u))\lambda(\gamma)f(\gamma').$$
Connes observes that the operator on the RHS is Hilbert-Schmidt as a
result of the convergence of the Poincare series:
$$\sum_{\gamma\in\Gamma} \exp(-2d(\gamma.u,x_0)).$$
Thus if
 $\gamma_0,\gamma_1,\gamma_2$ lie in $\Gamma$
then
$$\tr(\varepsilon\lambda(\gamma_0)[F,\lambda(\gamma_1)]
[F,\lambda(\gamma_2)])
=\frac{1}{2}{\tr}(\varepsilon[F,\lambda(\gamma_0)[F,\lambda(\gamma_1)]
[F,\lambda(\gamma_2)]).$$

Now $\lambda(\gamma_0)[F,\lambda(\gamma_1)]
[F,\lambda(\gamma_2)]$ is the operator
$$(\lambda(\gamma_0)[F,\lambda(\gamma_1)]
[F,\lambda(\gamma_2)]f)(\gamma)= \zeta(\gamma)\sigma(\gamma_1,\gamma_2)
\sigma(\gamma_0,\gamma_1\gamma_2)f((\gamma_0\gamma_1\gamma_2)^{-1}\gamma),$$
where  $\zeta(\gamma)$ denotes Clifford multiplication
by
$$(\varphi(\kappa_0^{-1}\gamma)-\varphi(\kappa_1^{-1}\gamma))(\varphi
(\kappa_1^{-1}\gamma)-\varphi(\kappa_2^{-1}\gamma)),
$$
with
$\kappa_j=\gamma_0\ldots\gamma_j$.
Finally we can now obtain a formula for the cyclic cocycle.
Following the calculation on page 344 of \cite{Co2}
we find that for $\gamma_0\gamma_1\gamma_2 \neq 1$ the
character of the  cocycle associated to our Fredholm module
is zero
while for $\gamma_0\gamma_1\gamma_2 =1$ it is given by
$$\tr(\varepsilon\lambda(\gamma_0)[F,\lambda(\gamma_1)]
[F,\lambda(\gamma_2)])= 2\sum_{\gamma\in\Gamma}
{\mbox{trace}}(\epsilon\zeta(\gamma))
{\sigma(\gamma_1,\gamma_2)},
$$
where `trace' denotes the matrix trace on the Clifford algebra
and we are utilising the fact that, for our choice of $\sigma$,
$$\sigma(\gamma_0,\gamma_1\gamma_2)=\sigma(\gamma_0,\gamma_0^{-1})=1.$$
Connes proves that
$\mbox{trace}(\varepsilon\zeta(\gamma))$
is the Euclidean
area of the triangle in the complex plane with
vertices corresponding to the tangent vectors
$\varphi(\kappa_j^{-1}\gamma)$.
Then the additive group cocycle on $\Gamma$ given by
$$c(1, \gamma_1,\gamma_1\gamma_2)= \sum_{\gamma\in
\Gamma}\mbox{trace}(\varepsilon\zeta(\gamma))$$
is what Connes calls the `volume' or area cocycle on $\Gamma$.
Thus we find that we have computed the character of our
Fredholm module to be:
$$\tau_c(\gamma_0,\gamma_1,\gamma_2)= c(1,\gamma_1,\gamma_1\gamma_2)
\sigma(\gamma_1,\gamma_2)$$
for $\gamma_0\gamma_1\gamma_2=1$, with $\tau_c$ being zero when
 $\gamma_0\gamma_1\gamma_2\neq 1$
(the normalisation differs from \cite{Co2}
page 295, but conforms with \cite{CM}). This formula extends to give
a non-trivial element of the cyclic cohomology
of the smooth subalgebra $A_{\sigma,g}$
via the formula
$$\tau_c(f^0,f^1,f^2)= \sum_{\gamma_0\gamma_1\gamma_2=1}
f^0(\gamma_0)f^1(\gamma_1)f^2(\gamma_2)c(1,\gamma_1,\gamma_1\gamma_2)\sigma(
\gamma_1,
\gamma_2),$$
for $f^0,f^1,f^2\in A_{\sigma,g}$.

Summarizing the discussion above, we have the first result of this section.

\begin{prop} There is a 2-summable Fredholm module
$(F,\mathcal{H})$ over
$A_{\sigma,g}$ whose Chern character is
given by the area cyclic 2-cocycle $\tau_c$. Therefore, by the
index pairing in \cite{Co2}, one has
\begin{gather*}
   \Index(P(F\otimes I)P) =    \langle[\tau_c], [P]\rangle , \\
   \end{gather*}
where $P$ denotes a projection in $A_{\sigma, g}\otimes \mathcal
K({\mathcal H}_1)$
and $\Index(P(F\otimes I)P)$ denotes the index of the Fredholm
operator $P(F\otimes I)P$.
\end{prop}

Assembling this proposition with our results from
Section 10 we have:

\begin{thm} Let $P$ denote a projection in $A_{\sigma, g}\otimes \mathcal
K({\mathcal H}_1)$
Then in the notation of Corollary 12 of the previous section, one has
$$
   \Index(P(F\otimes I)P) = 2(g-1)(\mbox{rank }\mathcal{E}^0 -
      \mbox{rank }\mathcal{E}^1)\in\mathbb{Z},
$$
where $\Index(P(F\otimes I)P)$ denotes the index of the Fredholm operator
$P(F\otimes I)P$ acting on the Hilbert space $P(\mathcal{H}\otimes
{\mathcal H}_1)$ and $\mu_\sigma^{-1}[P] = [{\npartial_{\mathcal{E}^0}^+}]
-
[{\npartial_{\mathcal{E}^1}^+}] \in K_0(\Sigma_g)$.
\end{thm}

\begin{cor} Let $P$ be a projection into a gap in the spectrum
of the discrete Hamiltonian $H_\tau$. Then $P\in A_{\sigma,g}$,
and may be regarded as a
twisted convolution operator by a function $p$ on $\Gamma$.
Then in the notation of Corollary 12:
$$
   \Index(PFP) =  \sum_{\gamma_0\gamma_1\gamma_2=1}
p(\gamma_0)p(\gamma_1)p(\gamma_2)c(1,\gamma_1,\gamma_1\gamma_2)\sigma(\gamma
_1,\gamma_2)
$$
$$= 2(g-1)(\mbox{rank }\mathcal{E}^0 -
      \mbox{rank }\mathcal{E}^1)\in\mathbb{Z}
$$
\end{cor}

Note that this explains the integrality of the cyclic 2-cocycle,
$$\sum_{\gamma_0\gamma_1\gamma_2=1}
p(\gamma_0)p(\gamma_1)p(\gamma_2)
c(1,\gamma_1,\gamma_1\gamma_2)\sigma(\gamma_1,\gamma_2),$$
in two different ways: firstly as the index of the Fredholm operator
$PFP$, and secondly as the topological index $2(g-1)(\mbox{rank}\mathcal{E}^0 -
      \mbox{rank}\mathcal{E}^1)$, which is also clearly an integer.

\section{The non-commutative unit disc}

In [Klim+Les1,2] Klimek and Lesznewski have introduced a non-commutative
unit disc
and higher genus Riemann surfaces.
Their disc algebra can be realised as a Toeplitz algebra obtained by
compressing the commutative algebra of functions on the disc using the
projection onto a holomorphic subspace of one of its representation spaces.
We shall describe their construction in a slightly more general setting.
The algebra $C_c(G/K)$ acts by multiplication ($f \mapsto M(f)$)
on $L^2(G/K,\mu)$ for any quasi-invariant measure $\mu$.
The group $G$ also has an induced $\sigma$ representation $W$ on this space,
and we shall suppose that there is an irreducible subrepresentation on a
subspace which is projected out by $P$.
(This is certainly true in the case considered in [Klim+Les1].)
The algebra $PM(C_c(G/K))P$ then gives the non-commutative analogue of
$C_c(G/K)$.
Now, by definition $G$ also acts and therefore defines automorphisms of this
algebra.
Since it commutes with $P$ the covariance algebra $PM(C_c(G/K))P\rtimes G$
is the same as $P(C_c(G/K)\rtimes G)P$, which is the compression of the
imprimitivity algebra $\ca = C_c(G/K)\rtimes G$.
For higher genus surfaces one simply takes the $\Gamma$-invariant part of
$PC_c(G/K)P$, which is consistent with our constructions above.

Suppose now that the irreducible subspace is defined by a reproducing kernel.
Invariance of the kernel means that it is defined by twisted convolution with
a continuous $\sigma$-positive definite function $\xi_P$ or, equivalently,
that
$$P= W(\xi_P) = \int \xi_P(g)W(g)\,dg.$$
Now observe that
 $\xi_P$ can be identified with an element of the imprimitivity
algebra so the covariance algebra can be identified with the compression
$\xi_P*\ca*\xi_P$ of the imprimitivity algebra.

In the cases of interest $\xi_P$ is the $\sigma$-positive-definite
function associated with a $C^\infty$-vector, and so is smooth.
This means that the natural module $\xi_P*\cm$ for
$\xi_P*\ca*\xi_P$ retains the structure of a Fredholm module.

\end{document}